\documentclass[12pt, reqno]{amsart}

\usepackage{mathtools}
\mathtoolsset{showonlyrefs}
\numberwithin{equation}{section}

\usepackage{graphicx}

\usepackage{amssymb}
\usepackage{amsthm}
\usepackage{amscd}
\usepackage{amsmath}
\usepackage{epic,eepic}
\usepackage{mathrsfs}
\usepackage{latexsym}
\usepackage{pifont}
\usepackage[all]{xy}
\usepackage {cancel}

\theoremstyle{plain}
\newtheorem{lemma}{Lemma}[section]
\newtheorem{prop}[lemma]{Proposition}
\newtheorem{thm}[lemma]{Theorem}
\newtheorem{cor}[lemma]{Corollary}

\theoremstyle{definition}

\newtheorem{rem}[lemma]{Remark}
\newtheorem{defi}[lemma]{Definition}
\newtheorem{exa}[lemma]{Example}

\newtheorem{problem}{Problem}

\newcommand{\bde}{\begin{defi}}
\newcommand{\ede}{\end{defi}\vspace{1mm}}
\newcommand{\ble}{\begin{lemma}}
\newcommand{\ele}{\end{lemma}}
\newcommand{\bpr}{\begin{prop}}
\newcommand{\epr}{\end{prop}}
\newcommand{\bt}{\begin{thm}}
\newcommand{\et}{\end{thm}}
\newcommand{\bco}{\begin{cor}}
\newcommand{\eco}{\end{cor}}
\newcommand{\bre}{\begin{rem}}
\newcommand{\ere}{\end{rem}}
\newcommand{\bex}{\begin{exa}}
\newcommand{\eex}{\end{exa}}
\newcommand{\bpf}{\begin{proof}}
\newcommand{\epf}{\end{proof}}

\newcommand{\mcC}{\mathcal{C}}
\newcommand{\mcD}{\mathcal{D}}
\newcommand{\mcE}{\mathcal{E}}
\newcommand{\mcF}{\mathcal{F}}
\newcommand{\mcG}{\mathcal{G}}
\newcommand{\mcH}{\mathcal{H}}

\newcommand{\mcM}{\mathcal{M}}

\newcommand{\mcO}{\mathcal{O}}

\newcommand{\mcS}{\mathcal{S}}
\newcommand{\mcT}{\mathcal{T}}

\newcommand{\mcV}{\mathcal{V}}

\newcommand{\mcZ}{\mathcal{Z}}

\newcommand{\mbC}{\mathbb{C}}

\newcommand{\mbF}{\mathbb{F}}

\newcommand{\mbM}{\mathbb{M}}

\newcommand{\mbP}{\mathbb{P}}
\newcommand{\mbQ}{\mathbb{Q}}
\newcommand{\mbR}{\mathbb{R}}
\newcommand{\mbS}{\mathbb{S}}

\newcommand{\mbV}{\mathbb{V}}

\newcommand{\mbZ}{\mathbb{Z}}

\newcommand{\mfS}{\mathfrak{S}}

\newcommand{\mfc}{\mathfrak{c}}

\newcommand{\mfl}{\mathfrak{l}}

\newcommand{\mfo}{\mathfrak{o}}

\newcommand{\mfs}{\mathfrak{s}}

\newcommand{\msE}{\mathscr{E}}

\newcommand{\msP}{\mathscr{P}}

\newcommand{\msU}{\mathscr{U}}

\newcommand{\msX}{\mathscr{X}}

\newcommand{\SSP}{\vspace{3mm}}
\newcommand{\LSP}{\vspace{5mm}}
\newcommand{\mr}{\mathrm}

\newcommand{\N}{N}

\newcommand{\LL}{[\![}
\newcommand{\RR}{]\!]}

\begin{document}

\title[Generalized hypergeometric equations and  dormant opers in char. $\leq 7$]{Generalized hypergeometric equations and \\ $2$d TQFT for dormant opers in characteristic  $\leq 7$}
\author{Keita Mori}
\author{Yasuhiro Wakabayashi}

\address{\emph{Keita Mori}
 \newline
 \textnormal{Graduate School of Information Science and Technology, The University of Osaka, Suita, Osaka 565-0871, JAPAN.
 }
 \newline
 \textnormal{\texttt{keita11262816@icloud.com}}}

 \address{\emph{Yasuhiro Wakabayashi}
 \newline
 \textnormal{Graduate School of Information Science and Technology, The University of Osaka, Suita, Osaka 565-0871, JAPAN.}
 \newline
  \textnormal{\texttt{wakabayashi@ist.osaka-u.ac.jp}}}

\date{}
\maketitle
\begin{abstract}
This note studies $\mathrm{PGL}_n$-opers arising from generalized hypergeometric differential equations in prime characteristic $p$. We prove  that these opers are rigid within the class of dormant opers. By combining this rigidity result with previous work in the enumerative geometry of dormant opers, we obtain a complete and explicit description of the $2$d TQFTs that compute  the number of  dormant $\mathrm{PGL}_n$-opers for primes $p \leq 7$.

\end{abstract}
\tableofcontents

\LSP
\section{Introduction} \label{S1}


A {\it $G$-oper} for a reductive group $G$ is a particular  type of    flat $G$-bundle on  an algebraic curve.
This notion was   introduced 
   in the context  of the geometric Langlands correspondence, serving as an element  for constructing Hecke eigensheaves  on the moduli space of bundles via quantization of  Hitchin's integrable system (cf. \cite{BeDr1}).
When $G = \mr{GL}_n$ or $\mr{PGL}_n$ (with $n \geq 2$),
such opers correspond to certain   
flat vector bundles of rank $n$ equipped with complete flags, and are associated with ordinary linear  differential operators  whose principal symbols are unit.
For instance,  each $\mr{PGL}_2$-oper on the projective line with  at most regular singularities at the three points $0$, $1$, and $\infty$  arises from a {\it Gauss hypergeometric differential operator} 
\begin{align} \label{Eq44}
D_{a, b, c} := \frac{d^2}{dx^2} + \left(\frac{c}{x} + \frac{1-c + a + b}{x-1} \right) \cdot \frac{d}{dx} + \frac{ab}{x (x-1)},
\end{align}
determined by a triple  $(a, b, c)$ of parameters, where $x$ denotes the standard  coordinate on the projective line.

In the complex analytic  setting, 
$\mr{PGL}_2$-opers on a closed   Riemann  surface
can be identified, via the Riemann-Hilbert correspondence,  with certain refinements of its complex structure known as   
 {\it projective structures}, i.e.,  atlases  of coordinate charts  whose  transition functions are   M\"{o}bius transformations.
A canonical example of a projective structure is constructed by the system of local inverse maps of the universal covering arising from uniformization.

 $G$-opers  {\it in prime characteristic $p > 0$} have been studied as a part of   the  characteristic-$p$ analogue  of the geometric Langlands correspondence (cf. ~\cite{BeTr}), as well as  in connection with various topics, including 
    $p$-adic Teichm\"{u}ller theory (cf., e.g.,  ~\cite{Moc1},  \cite{Moc2},  ~\cite{JRXY}, ~\cite{JoPa}, ~\cite{LaPa}, and  ~\cite{LiOs}).
 A central concept in   these developments is the  {\it $p$-curvature} of a flat $G$-bundle, which serves as an  
 invariant measuring  the obstruction to compatibility between  $p$-power operations  on certain spaces of   infinitesimal symmetries.
 This invariant also involves 
the Grothendieck-Katz conjecture, which provides  a conjectural  criterion for  
 the algebraicity of solutions to linear differential equations 
  (cf. ~\cite{NKa3}, ~\cite{And}).
  
A $G$-oper is said to be  {\it dormant} if its $p$-curvature  vanishes.
In the context of $p$-adic Teichmuller theory, dormant $\mr{PGL}_2$-opers (or more generally, $\mr{PGL}_2$-opers with nilpotent $p$-curvature) may  be viewed  as analogues of ``well-behaved" projective structures on Riemann surfaces such as those arising from uniformization.
The theory  of dormant $G$-opers for general $G$ has been   developed extensively  in the author's works (cf., e.g.,  ~\cite{Wak1}, ~\cite{Wak2}, ~\cite{Wak3}, ~\cite{Wak5}, ~\cite{Wak7}, \ ~\cite{Wak8}, ~\cite{Wak9},  and  ~\cite{Wak10}).

Now, let us consider the case where   $G = \mr{PGL}_n$ with $2 \leq n \leq p$, and fix a pair of nonnegative integers $(g, r)$ satisfying  $2g-2 +r >0$.
A  central object   in aforementioned works is the moduli stack 
\begin{align} \label{Eq223}
\mcO p^{^\mr{Zzz...}}_{n, \rho, g, r}
\end{align}
(cf. \eqref{Eq100}), which    classifies  pairs $(\msX, \msE^\spadesuit)$
 consisting of a pointed curve $\msX$ in $\overline{\mcM}_{g, r}$ (:= the moduli stack of $r$-pointed stable curves of genus $g$ in characteristic $p$) and 
 a dormant $\mr{PGL}_n$-oper $\msE^\spadesuit$ on $\msX$ of prescribed radii $\rho$ (cf. ~\cite[Definition 2.32]{Wak5}  for the definition of radius).
As shown  in ~\cite{Wak2} and ~\cite{Wak5},
the stack $\mcO p^{^\mr{Zzz...}}_{n,  \rho, g, r}$ is finite and generically \'{e}tale over $\overline{\mcM}_{g, r}$, so  it is meaningful to consider  its  generic degree
 \begin{align}
 N_{p, n, \rho, g, r}  := \mr{deg}(\mcO p^{^\mr{Zzz...}}_{n, \rho, g, r}/\overline{\mcM}_{g, r})
 \end{align}
 (cf. \eqref{Eq113}),
  which  counts  the number of dormant $\mr{PGL}_n$-opers of radii $\rho$ on a general curve in $\overline{\mcM}_{g, r}$.

Recall from   ~\cite{Wak20} (or ~\cite{Wak11}) that 
 the values $N_{p, n, \rho, g, r}$
 satisfy a factorization rule governed by various 
  gluing  procedures of underlying stable  curves,
  and  this structure endows them with the properties of   a $2$-dimensional topological quantum field theory ($2$d TQFT).
A major goal of our  study is to understand this $2$d TQFT, as it provides a bridge between
 the theory of dormant opers and other enumerative geometries, such as  the Gromov-Witten theory of Grassmannians and the conformal field theory associated to  the affine Lie algebras (cf. ~\cite{Wak5}).

Concerning this, 
it is known that a $\mr{GL}_n$-oper or a $\mr{PGL}_n$-oper  is dormant if and only if the corresponding  differential operator  {\it admits  a full set of solutions}.
When a given $\mr{PGL}_2$-oper arises from the operator $D_{a, b, c}$ as defined  above,
this condition translates  into the requirement   that
   $(a, b, c)$ is  the mod $p$ reduction of a triple of integers $(\widetilde{a}, \widetilde{b}, \widetilde{c})$ in $\{1, \cdots, p \}$ satisfying  either $\widetilde{a} < \widetilde{c} \leq \widetilde{b}$ or $\widetilde{b} < \widetilde{c} \leq \widetilde{a}$ (cf. ~\cite{Iha},  ~\cite{Kat2}).
This characterization  enables   an explicit description of the $2$d TQFT governing dormant $\mr{PGL}_2$-opers.

However, a comprehensive understanding of this TQFT remains out of reach for general $n$, 
as little seems to be known beyond the special cases $n= p-2, p-1, p$.
Thus, in this note, we take a step toward a broader understanding by  investigating the number of dormant $\mr{PGL}_n$-opers for $n >2$.

Thanks to the factorization property of the values $N_{p, n, \rho, g, r}$,
 it suffices to consider the case $(g, r) = (0, 3)$.
 A key insight in our discussion is that, in this case, most dormant opers arise from {\it generalized hypergeometric differential operators}, expressed as
\begin{align}
D_{\alpha, \beta} := x \frac{d}{dx} \cdot  \prod_{j=1}^{n-1} \left(x \frac{d}{dx}  + \beta_j -1\right) - x \cdot  \prod_{j=1}^n \left(x \frac{d}{dx}  + \alpha_j\right)
\end{align}
for suitable parameters   $\alpha := (\alpha_1, \cdots, \alpha_n)$, $\beta := (\beta_1, \cdots, \beta_{n-1})$.
We then  apply a result of Katz (cf.   ~\cite{Kat5}), which  generalizes the case $n=2$, i.e,   Gauss' hypergeometric operators, to determine the conditions under which such an operator  $D_{\alpha, \beta}$ has a full set of root functions (i.e.,  functions annihilated by $D_{\alpha, \beta}$), or equivalently, when the associated $\mr{PGL}_n$-oper is dormant.
As a consequence, we obtain a complete and explicit description of   the $2$d TQFT for dormant $\mr{PGL}_n$-opers  in characteristic $p \leq 7$ (cf. Section \ref{SS7} for details).
This provides  the first effective method for computing the values $N_{p, n, \rho, g, r}$ in the previously unexplored range   $2 < n < p-2$  and $r > 0$.

\LSP
\subsection*{Notation and Conventions} \label{SS0}

Throughout  this paper, we fix 
an odd  prime number  $p$, an algebraically closed  field $k$ of characteristic $p$, and  an integer $n$ 
  with $1 < n < p$.
We denote by $\mr{GL}_n$ (resp., $\mr{PGL}_n$) the general (resp.,  projective) linear group of $k^{n}$.

For a vector bundle (i.e., a locally free coherent sheaf)  $\mcF$ on a scheme $S$, we denote by $\mbV (\mcF)$ the relative affine scheme over $S$  associated to $\mcF$, i.e., the spectrum 
\begin{align}
\mbV (\mcF) := \mcS pec (\mr{Sym}_{\mcO_S} (\mcV^\vee)),
\end{align}
where 
  $\mr{Sym}_{\mcO_S}(\mcF^\vee)$ denotes the symmetric algebra  of $\mcF^\vee$ over $\mcO_S$.

If $S$ is a scheme over  $k$,
then we denote by $S^{(1)}$ its  Frobenius twist over $k$, i.e., 
the base-change of $S$ along the absolute Frobenius endomorphism  of $\mr{Spec}(k)$.
Let  $F_{S/k} : S \rightarrow S^{(1)}$ denote the relative Frobenius morphism of $S/k$.

\vspace{10mm}
\section{$\mr{PGL}_n$-opers induced  from generalized hypergeometric operators} \label{S9}

In this section, we study  generalized hypergeometric differential operators  in characteristic $p$, as well as the  
$\mr{PGL}_n$-opers 
induced from them.
These opers are described in terms of logarithmic connections on vector bundles, following the approach of   ~\cite{Wak5}).
We emphasize that, on the $3$-pointed projective line,  the dormant opers, i.e., those with vanishing $p$-curvature, are 
classified by  certain configurations of radii.
This classification  follows from 
Katz's result characterizing when a generalized hypergeometric operator admits  a full set of root functions  (cf. Proposition \ref{Prop4}).
The main result of this section asserts  that  such dormant  $\mr{PGL}_n$-opers are uniquely  determined by their radii  (cf. Theorem \ref{Prop90}).

\LSP
\subsection{Generalized hypergeometric operators in characteristic $p$} \label{SS1}

Let  $k (x)$ denote  the field of rational functions in the variable $x$  over $k$, and endow it with the structure of a differential field over $k$ via the derivation
 $\delta_x := x \frac{d}{dx}$.

Let $n$, $m$ be positive integers and consider tuples   $\alpha := (\alpha_1, \cdots, \alpha_n) \in k^{n}$, $\beta := (\beta_1, \cdots, \beta_m) \in k^{m}$.
To this pair $(\alpha, \beta)$, we associate 
the   {\bf generalized hypergeometric differential operator} 
\begin{align} \label{Eq6}
D_{\alpha, \beta} := \delta_x \cdot  \prod_{j=1}^{m} \left(\delta_x + \beta_j -1 \right) - x \cdot \prod_{j=1}^n \left(\delta_x + \alpha_j\right),
\end{align}
defined as a linear differential operator on $k (x)$.
When we regard $k (x)$ as a $k (x^p)$-vector space with basis $1, x, \cdots, x^{p-1}$,
the operator  $D_{\alpha, \beta}$ defines a  $k(x^p)$-linear endomorphism of $k (x)$.
In particular, the kernel $\mr{Ker}(D_{\alpha, \beta})$ forms  a $k(x^p)$-vector subspace of $k (x)$, and its dimension  satisfies  $\mr{dim}_{k (x^p)} (\mr{Ker}(D_{\alpha, \beta})) \leq p$.
In what follows,  we investigate  how  the dimension of  $\mr{Ker}(D_{\alpha, \beta})$ can be described in terms of the data  $(\alpha, \beta)$.

Note that the operator $D_{\alpha, \beta}$ is invariant under reordering of the entries in  $\alpha$ and $\beta$.
Thus, without loss of generality,  we may assume the following conditions after possibly reordering the elements:
\begin{itemize}
\item
$\alpha_1, \cdots, \alpha_{n'}, \beta_1, \cdots, \beta_{m'} \in \mbF_p$ and $\alpha_{n' +1}, \cdots, \alpha_n, \beta_{m'+1}, \cdots, \beta_{m} \in  k \setminus \mbF_p$, for some integers   $0 \leq n' \leq n$ and $0 \leq m' \leq m$; 
\item
The following  inequalities are fulfilled:
 \begin{align} \label{Eq134}
 p \geq \widetilde{\alpha}_1 \geq \widetilde{\alpha}_2 \geq \cdots \geq \widetilde{\alpha}_{n'} \geq 1 \ \ \  \text{and}  \ \ \ p \geq \widetilde{\beta}_1 \geq \widetilde{\beta}_2 \geq \cdots \geq \widetilde{\beta}_{m'} \geq 1,
 \end{align}
where, for each $\gamma \in \mbF_p$, we denote by  $\widetilde{\gamma}$ the unique integer in  $\{ 1, \cdots, p \}$ congruent to $\gamma$ modulo $p$.
\end{itemize}
We then define 
\begin{align} \label{Eq401}
T_{\alpha, \beta}
\end{align}
to be the subset of $\{1, \cdots, m' \}$ consisting of those indices  $j$ for which  there exists some $j'$ satisfying  $\widetilde{\beta}_j > \widetilde{\alpha}_{j'} \geq \widetilde{\beta}_{j+1}$, where we set $\widetilde{\beta}_{m' +1} := 1$ by convention.

Now, observe that, for each $s \in \mbZ_{\geq 0}$, the  following equality holds:
\begin{align}
\left(\frac{1}{x} \cdot D_{\alpha, \beta} \right) (x^s) =  \left(s \cdot \prod_{j=1}^{m} (s-1+ \beta_j)\right) x^{s-1} + \left( -\prod_{j=1}^n (s +\alpha_j) \right)x^{s}.
\end{align}
Define the polynomials 
\begin{align}
P (X) := -\prod_{j=1}^n (X + \alpha_j), \hspace{10mm} Q (X) := (X+1) \cdot \prod_{j=1}^{m} (X + \beta_j).
\end{align}
Then,   the matrix representation of the operator $\frac{1}{x} \cdot D_{\alpha, \beta}$ with respect to the basis $1, x, \cdots, x^{p-1}$ of the $k (x^p)$-vector space $k (x)$ is given by the upper bidiagonal matrix
\begin{align}
R_{\alpha, \beta} := \begin{pmatrix}
P (0)& Q (0) & 0 & 0 & 0& \cdots & 0
\\
0 & P (1) & Q (1)& 0& 0 & \cdots& 0
\\
0& 0& P (2)& Q (2)& 0& \cdots & 0
\\
\vdots & \vdots && \ddots &\ddots &\ddots & \vdots 
\\
\vdots & \vdots   &&&P (p-3)& Q (p-3)& 0
\\
0 & 0   & \cdots & \cdots &0 &P (p-2)& Q (p-2)
\\
0 & 0 & \cdots & \cdots & 0 & 0 & P (p-1)
 \end{pmatrix}.
\end{align}
In particular, we have the identity 
\begin{align} \label{Eq139}
p- \mr{rank} \left(R_{\alpha, \beta}\right)  =  \mr{dim}_{k (x^p)} \left(\mr{Ker} \left(\frac{1}{x} \cdot D_{\alpha, \beta}\right)\right)  =  \mr{dim}_{k (x^p)} (\mr{Ker} (D_{\alpha, \beta})).
\end{align}
Moreover, note that the diagonal entry $P (\ell)$ (for $0 \leq \ell \leq p-1$) is congruent to zero modulo $p$ precisely when 
   $\ell = p- \widetilde{\alpha}_j$ for some $j \in \{1, \cdots, n' \}$.

  \SSP
\bpr \label{Prop2}
For $(\alpha, \beta) \in k^{n} \times k^{m}$ as above,  the following equality holds:
\begin{align}
\mr{rank} (\mr{Ker}(D_{\alpha, \beta})) = \sharp (T_{\alpha, \beta}).
\end{align}
\epr
\begin{proof}
For convenience, 
we set   $b_{0} := -1$, $b_{m'+1} := p-1$, and $b_j := p- \widetilde{\beta}_j$ for $j=1, \cdots, m'$.
For each $j= 1, \cdots, m'+1$, we define a  $k (x^p)$-vector subspace of $k (x)$  by 
 $L_j := \bigoplus_{s = b_{j-1}+1}^{b_{j}}k (x^p) x^{s} \left(\subseteq L \right)$.
Then, $k (x)$ decomposes as  $L = \bigoplus_{j=1}^{m'+1}L_j$, which induces  a decomposition $\frac{1}{x} \cdot D_{\alpha, \beta} = \bigoplus_{j=1}^{m' +1} D'_j$,
where each $D'_j$ is  a $k(x^p)$-linear endomorphism of $L_j$.
The matrix representing  $D_j$ with respect to the basis  $x^{b_{j-1}+1}, \cdots, x^{b_j}$ is given by  the $(b_j - b_{j-1}) \times (b_j - b_{j-1})$  bidiagonal matrix
\begin{align}
R'_j := \begin{pmatrix}
P (b_{j-1}+1)& Q (b_{j-1}+1) & 0 & \cdots & 0& 0 & 0
\\
0 & P (b_{j-1}+2) & Q (b_{j-1}+2)& \cdots & 0 & 0 & 0
\\
0& 0& P (b_{j-1}+3)& \ddots & 0 & 0 & 0
\\
\vdots & \vdots && \ddots &\ddots &\ddots & \vdots 
\\
0& 0   &0 &\cdots &P (b_{j}-2)& Q (b_{j}-2)& 0
\\
0 & 0   & 0 & \cdots &0 &P (b_{j}-1)& Q (b_{j}-1)
\\
0 & 0 & 0 & \cdots & 0 & 0& P (b_j) 
 \end{pmatrix}.
\end{align}
Since all off-diagonal entries  $Q (b_{j-1}+1), \cdots, Q (b_j -1)$  are nonzero,  we have $\mr{rank} (R'_j) = b_j - b_{j-1} -1$ if  $j-1 \in T_{\alpha, \beta}$ ($\Leftrightarrow b_{j-1} < p-\widetilde{\alpha}_{j'} \leq b_j$ for some $j'$), and $\mr{rank} (R'_j) = b_j - b_{j-1}$ if otherwise.
Hence, the following equalities hold:
\begin{align}
\mr{dim}_{k (x^p)} \left(\mr{Ker}\left(\frac{1}{x} \cdot D_{\alpha, \beta}\right)\right) & = p- \mr{rank} R_{a, b} \\
& =p-  \sum_{j=1}^{m' +1} \mr{rank} (R'_j) \notag \\
& = p- \left( \left( \sum_{j=1}^{m' +1} (b_j - b_{j-1}) \right) - \sharp (T_{\alpha, \beta})\right) \notag \\
& = p - \left(b_{m'} - b_{0} - \sharp (T_{\alpha, \beta}) \right) \notag \\
&=   \sharp (T_{\alpha, \beta}).
  \notag
\end{align}
The assertion then follows from \eqref{Eq139}.
\end{proof}
\SSP

The following assertion is  a direct consequence of the  proposition above.
It was already established  in 
 ~\cite[Sublemma 5.5.2.1]{Kat5}  (and ~\cite[Section 1.6]{Iha} for  the case of Gauss' hypergeometric operators) to provide a complete classification of hypergeometric differential operators with finite monodromy; see also   ~\cite[Remark 4.9]{BeHe}. 

\SSP
\bco \label{Cor11}
(Recall that $(\alpha, \beta)$ has assumed to satisfy the inequalities in \eqref{Eq134}.)
The kernel $\mr{Ker}(D_{\alpha, \beta})$ has rank $n$ (as a $k (x^p)$-vector space) if and only if the following two conditions are fulfilled:
\begin{itemize}
\item[(1)]
$n' =n$ and $m' = m$, i.e., the two sets $\{ \alpha_j \}_{j=1}^n, \{ \beta_j \}_{j=1}^m$ are contained in $\mbF_p$;
\item[(2)]
$m = n-1$ and the following chain of inequalities  holds:
\begin{align}
\widetilde{\alpha}_1 \geq \widetilde{\beta}_1 > \widetilde{\alpha}_2 \geq \widetilde{\beta}_2 > \cdots \geq \widetilde{\beta}_{n-1}> \widetilde{\alpha}_n.
\end{align}
\end{itemize}
\eco
\begin{proof}
By the definition of $T_{\alpha, \beta}$,
the equality $\sharp (T_{\alpha, \beta}) = n$ holds precisely when
both conditions  (1) and (2) are fulfilled.
Thus, the assertion follows from  Proposition \ref{Prop2}.
\end{proof}

\LSP
\subsection{Dormant $\mr{PGL}_n$-opers} \label{SS2}

We now begin our discussion of  
dormant $\mr{PGL}_n$-opers (= dormant $\mfs \mfl_n$-opers) of prescribed radii on pointed  curves.
To simplify the exposition, we will work with  equivalent objects described in terms of log connections on vector bundles (without using the formulation by log structures).
For a comprehensive treatment of  $\mr{PGL}_n$-opers  on log curves and their various properties, we refer the reader to ~\cite{Wak5}.

Let $(g, r)$ be a pair of nonnegative integers with $2g-2+r >0$, and
 $\msX := (X, \{ \sigma_i \}_{i=1}^r)$  an $r$-pointed proper  smooth curve of genus $g$ over $k$, where $X$ denotes the underlying  curve and $\sigma_1, \cdots, \sigma_r$ are ordered  marked points on $X$.
Denote by $\Omega \left(:= \Omega_{X/k} (\sum_{i=1}^r \sigma_i) \right)$ the sheaf of logarithmic $1$-forms on $X/k$ with poles along the marked points $\sigma_i$.
Also, denote by $\mcT$ its dual, i.e., $\mcT := \Omega^\vee$.
For each $j \in \mbZ_{\geq 0}$, we have  
 the sheaf of crystalline logarithmic differential operators $\mcD_{\leq j}$ on $X/k$ (with poles along $\sigma_i$'s,  as above) of order $\leq j$ (cf.  ~\cite[D\'{e}finition 2.3.1]{Mon}, ~\cite[Section 4.2.1]{Wak5}).
We  set  $\mcD := \bigcup_{j \in \mbZ_{\geq 0}} \mcD_{\leq j}$.

Recall that  a {\bf log connection} on an $\mcO_X$-module $\mcF$ is a $k$-linear morphism $\nabla : \mcF \rightarrow \Omega \otimes \mcF$ satisfying the usual Leibnitz rule
\begin{align}
\nabla_\partial (a \cdot v) = \partial (a) \cdot v + a \cdot \nabla_\partial (v)
\end{align}
 for any local sections $\partial \in \mcT$, $a \in \mcO_X$,  and $v \in \mcF$, where $\nabla_{\partial} := (\partial \otimes \mr{id}_\mcF) \circ \nabla$ (cf. ~\cite[Definition 4.1]{Wak5}).

The {\bf $p$-curvature} of such a log connection $\nabla$ is defined as the $\mcO_X$-linear morphism
\begin{align}
\psi (\nabla) : \mcT^{\otimes p} \rightarrow \mcE nd_{\mcO_X} (\mcF)
\end{align}
 determined uniquely by the condition that $\psi (\nabla) (\partial^{\otimes p}) = \nabla_{\partial}^p - \nabla_{\partial^{p}}$ for any local section $\partial \in \mcT$, where $\partial^{p}$ denotes
  the local section of $\mcT$ corresponding to the $p$-th iterate of the (locally defined) derivation on $\mcO_X$ associated to  $\partial$.

We now  fix 
an {\bf  $n$-theta characteristic} of $\msX$ in the sense of ~\cite[Definition 4.31, (i)]{Wak5}), i.e.,  a pair 
\begin{align} \label{Eq404}
\vartheta := (\varTheta, \nabla_\vartheta)
\end{align}
consisting of a line bundle $\varTheta$ on $X$ and a log connection $\nabla_\vartheta$ on  the line bundle $\mcT^{\otimes \frac{n (n-1)}{2}} \otimes \varTheta^{\otimes n}$.
Moreover, we assume that $\nabla_\vartheta$ has vanishing $p$-curvature.
 (Such an $n$-theta characteristic always exists, according to the discussion in ~\cite[Section 4.6.4]{Wak5}.)
 Then,
the residue $\mr{Res}_{\sigma_i} (\nabla_\vartheta)$ of $\nabla_\vartheta$ at each marked point $\sigma_i$ ($i=1, \cdots, r$) is given by an element of $\mbF_p$ (cf., e.g.,  ~\cite[Proposition-Definition 4.8]{Wak20}).

We define   
\begin{align}
\mcF_\varTheta := \mcD_{\leq n-1} \otimes \varTheta \ \ \  \text{and}  \ \ \ \mcF_\varTheta^j := \mcD_{\leq n - j-1} \otimes \varTheta \  (j=0, \cdots, n).
\end{align} 
In particular, $\{ \mcF_\varTheta^j \}_{j=0}^n$ forms  an $n$-step decreasing filtration on $\mcF_\varTheta$ whose graded pieces  are line bundles.
The determinant $\mr{det}(\mcF_\varTheta)$ of $\mcF_\varTheta$ admits  a sequence of canonical isomorphisms
\begin{align}
\label{Eq22}
\mr{det}(\mcF_\varTheta) \xrightarrow{\sim} \bigotimes_{j=0}^{n-1} \mcF^j_\varTheta / \mcF^{j+1}_\varTheta \xrightarrow{\sim} \bigotimes_{j=0}^{n-1} \left( \mcT^{\otimes n-j-1} \otimes \varTheta \right) \xrightarrow{\sim} \mcT^{\otimes \frac{n(n-1)}{2}} \otimes \varTheta^{\otimes n}.
\end{align}

\SSP
\bde[cf. ~\cite{Wak5}, Definition 4.36] \label{Def3}
\begin{itemize}
\item[(i)]
 A {\bf  $(\mr{GL}_n, \vartheta)$-oper} on $\msX$ is a log connection $\nabla^\diamondsuit$ on $\mcF_\varTheta$ satisfying  the following three conditions:
\begin{itemize}
\item
For each $j=1, \cdots, n-1$, $\nabla^\diamondsuit (\mcF^j_\varTheta)$ is contained in $\Omega \otimes \mcF^{j-1}_\varTheta$;
\item
For each $j=1, \cdots, n-1$, the well-defined $\mcO_X$-linear morphism
\begin{align}
\mr{KS}^j : \mcF^j_\varTheta/\mcF^{j+1}_\varTheta \rightarrow \Omega \otimes (\mcF^{j-1}_\varTheta/\mcF^j_\varTheta)
\end{align}
given  by $\overline{a} \mapsto \overline{\nabla^\diamondsuit (a)}$ for any local section $a \in \mcF^j_\varTheta$ (where $\overline{(-)}$'s denote the images in the respective quotients) is an isomorphism;
\item
The log connection $\mr{det}(\nabla^\diamondsuit)$ on $\mr{det}(\mcF_\varTheta)$  induced by $\nabla^\diamondsuit$ commutes  with $\nabla_\vartheta$ via \eqref{Eq22}.
\end{itemize}

\item[(ii)]
A $(\mr{GL}_n, \vartheta)$-oper $\nabla^\diamondsuit$ is said to be {\bf dormant} if its $p$-curvature vanishes.

\item[(iii)]
Two (dormant) $(\mr{GL}_n, \vartheta)$-opers $\nabla_\circ^\diamondsuit$, $\nabla_\bullet^\diamondsuit$ are said to be  {\bf isomorphic}
if there exists an $\mcO_X$-linear  automorphism of $\mcF_\varTheta$ preserving the filtration $\{ \mcF_\varTheta^j \}_j$ such that 
$\nabla_\circ^\diamondsuit$ commutes with $\nabla_\bullet^\diamondsuit$ via this automorphism.
\end{itemize}
\ede
\SSP

Next, for $R \in \{ \mbF_p, k \}$, we denote  by $\Delta_R$ the image of the diagonal embedding $R \hookrightarrow R^{n}$, which is a group homomorphism.
In particular, this yields  the quotient set $R^{n}/\Delta_R$.
The symmetric group of $n$ letters $\mfS_n$
acts on 
the set $R^{n}$ by permuting the entries of each tupes.
This action descends to a  well-defined $\mfS_n$-action on $R^{n}/\Delta_R$.
Accordingly, we obtain the quotient sets $\mfS_n \backslash R^{n}$ and
 \begin{align} \label{Eq403}
 \mfc (R) := \mfS_n \backslash R^{n}/\Delta_R.
 \end{align}
 Note that each element of  $\mfS_n \backslash R^{n}$ can be interpreted as a multiset of elements of $R$ with cardinality $n$.
 The natural projection $R^n \twoheadrightarrow R^{n}/\Delta_R$ induces a surjection $\mfS_n \backslash R^{n} \twoheadrightarrow \mfc (R)$.

Given  elements $s_1, \cdots, s_n \in R$,
we write $\LL s_1, \cdots, s_n \RR$ for the element of   $\mfc (R)$ represented by the $n$-tuple $s := (s_1, \cdots, s_n)$.
For each element $\overline{s} := \LL s_1, \cdots, s_n \RR$ (where $s_1, \cdots, s_n \in k$) of 
 $\mfc (k)$,
   the  diagonal $n \times n$ matrix with diagonal entries $s_1, \cdots, s_n$
    specifies a well-defined  element $\rho_{\overline{s}}$ in  
 the GIT quotient  of $\mr{Lie}(\mr{PGL}_n)$ by the adjoint $\mr{PGL}_n$-action.
 The resulting assignment $\overline{s} \mapsto \rho_{\overline{s}}$ yields  a natural   identification between    $\mfc (k)$ and  this GIT quotient.
By the inclusion $\mbF_p \hookrightarrow  k$,  we regard  $\mfc(\mbF_p)$ as a subset of $\mfc (k)$.

Define   $\widetilde{\Xi}_{p, n}$ to be the subset of $\mfS_n \backslash \mbF^{r}_p$ consisting of multisets $[d_1, \cdots, d_n]$ in which  the elements $d_1, \cdots, d_n$ are pairwise  distinct.
(In particular, $\widetilde{\Xi}_{p, n}$ can be identified with the set  of all $n$-element subsets of $\mbF_p$.)
We denote the image of $\widetilde{\Xi}_{p, n}$ via the projection $\mfS_n \backslash \mbF^{n}_p \twoheadrightarrow \mfc (\mbF_p)$ by 
\begin{align} \label{Eq410}
\Xi_{p, n}.
\end{align}

Let us take an $r$-tuple $\rho := (\rho_i)_{i=1}^r$ of elements of 
$\mfc (k)$.
Since  $n < p$,
 each $\rho_i$ is  represented by a unique  $n$-tuple  $(a_{i, 1}, \cdots, a_{i, n}) \in k^{n}$ such that the sum $\sum_{i=1}^n a_{i, n}$ coincides with  $\mr{Res}_{\sigma_i} (\nabla_\vartheta)$.
 If $\rho_i \in \mfc (\mbF_p)$, then  this $n$-tuple belongs to $\mbF_p^{n}$.

\SSP
\bde[cf. ~\cite{Wak5}, Definition 4.43] \label{Def6}
We say that a  $(\mr{GL}_n, \vartheta)$-oper $\nabla^\diamondsuit$ is {\bf  of radii $\rho$}
if, for each $i=1, \cdots, r$, the characteristic polynomial of the residue matrix  of $\nabla^\diamondsuit$ at $\sigma_i$ coincides with that of  the diagonal matrix with diagonal entries $a_{i, 1}, \cdots, a_{i, n}$.
(When  $r = 0$, any $(\mr{GL}_n, \vartheta)$-oper is said to be  of radii $\emptyset$.)
\ede
\SSP

For a $(\mr{GL}_n, \vartheta)$-oper $\nabla^\diamondsuit$ of radii $\rho$ on $\msX$,
its projectivization (in a certain natural sense) determines a  $\mr{PGL}_n$-oper $\nabla^{\diamondsuit \Rightarrow \spadesuit}$ of radii $\rho$ on $\msX$
 (cf. ~\cite[Definitions 2.1, 4.43]{Wak5} for the definition of a $\mr{PGL}_n$-oper of prescribed radii);
the resulting assignment $\nabla^\diamondsuit \mapsto  \nabla^{\diamondsuit \Rightarrow \spadesuit}$ gives a bijective correspondence 
\begin{align} \label{Eq14}
\begin{pmatrix}\text{the set of isomorphism} \\  \text{classes of $(\mr{GL}_n, \vartheta)$-opers} \\\text{(resp., dormant $(\mr{GL}_n, \vartheta)$-opers)} \\ \text{ of radii $\rho$ on $\msX$} \end{pmatrix}
\xrightarrow{\sim}
\begin{pmatrix}\text{the set of isomorphism} \\  \text{classes of $\mr{PGL}_n$-opers} \\\text{(resp., dormant $\mr{PGL}_n$-opers)} \\ \text{ of radii $\rho$ on $\msX$} \end{pmatrix}
\end{align}
(cf. ~\cite[Theorem 4.66, Corollary 4.70]{Wak5}).

It follows from ~\cite[Proposition 6.14]{Wak20}
that each element of $\mfc (k)$ arising  as the radius of 
 a dormant $\mr{PGL}_n$-oper must  lie in $\Xi_{p, n}$.

\LSP
\subsection{Correspondence with differential operators} \label{SS18}

We define 
\begin{align}
\mcD {\it iff}_{\!\vartheta, \leq n} := \mcH om_{\mcO_X} (\varTheta^\vee, (\Omega^{\otimes n} \otimes \varTheta^\vee) \otimes \mcD_{\leq n}).
\end{align}
As discussed in ~\cite[Remark 4.2]{Wak5},
each global section of this sheaf can be regarded as an $n$-th order  linear differential operator from $\varTheta^\vee$ to $\Omega^{\otimes n} \otimes \varTheta^\vee$.
This sheaf admits   the composite surjection 
\begin{align} \label{Eq20}
\Sigma : \mcD {\it iff}_{\!\vartheta, \leq n} \twoheadrightarrow  \mcH om_{\mcO_X} (\varTheta^\vee, \Omega^{\otimes n} \otimes \varTheta^\vee \otimes \mcT^{\otimes n}) \xrightarrow{\sim} \mcO_X,
\end{align}
where  the first and second arrows  arise from the canonical isomorphisms $\mcD_{\leq n}/\mcD_{\leq n-1} \xrightarrow{\sim} \mcT^{\otimes n}$ and $\Omega^{\otimes n} \otimes \mcT^{\otimes n} \xrightarrow{\sim} \mcO_X$, respectively.
For each global section  $D \in H^0 (X, \mcD {\it iff}_{\!\vartheta, \leq n})$, we refer to $\Sigma (D) \in H^0 (X, \mcO_X) \left(= k \right)$ as the {\bf principal symbol} of $D$.

Let us take a global section $D^\clubsuit : \varTheta^\vee \rightarrow (\Omega^{\otimes n} \otimes \varTheta^\vee) \otimes \mcD_{\leq n}$ of the inverse image  $\Sigma^{-1} (1)$.
This  corresponds to an $\mcO_X$-linear morphism $D' : \mcT^{\otimes n} \otimes \varTheta \rightarrow \mcD_{\leq n} \otimes \varTheta$ via 
 the composite of natural isomorphisms 
\begin{align}
\mcD {\it iff}_{\!\vartheta, \leq n} \xrightarrow{\sim} \Omega^{\otimes n} \otimes \varTheta^\vee \otimes \mcD_{\leq n} \otimes \varTheta \xrightarrow{\sim} \mcH om_{\mcO_X} (\mcT^{\otimes n} \otimes \varTheta, \mcD_{\leq n} \otimes \varTheta).
\end{align} 
Using $D'$, we construct 
 the left $\mcD$-module $(\mcD \otimes\varTheta) / \langle \mr{Im}(D') \rangle$, i.e., the quotient of the $\mcD$-module $\mcD \otimes \varTheta$ by the $\mcD$-submodule generated by the sections of $\mr{Im}(D')$.
Since $\Sigma (D^\clubsuit) =1$,
the composite
\begin{align}
\mcF_\varTheta \xrightarrow{\mr{inclusion}} \mcD \otimes \varTheta \xrightarrow{\mr{quotient}} (\mcD \otimes\varTheta) / \langle \mr{Im}(D') \rangle
\end{align}
is an isomorphism of $\mcO_X$-modules.
The $\mcD$-action on $(\mcD \otimes\varTheta) / \langle \mr{Im}(D') \rangle$ determines, via this composite isomorphism, a log connection
\begin{align}
D^{\clubsuit \Rightarrow \diamondsuit} : \mcF_\varTheta \rightarrow \Omega \otimes \mcF_\varTheta
\end{align}
on $\mcF_\varTheta$ (cf. ~\cite[Section 4.2.2]{Wak5}).
It is immediately verified  that $D^{\clubsuit \Rightarrow \diamondsuit}$ defines a $(\mr{GL}_n, \varTheta)$-oper, in the sense of ~\cite[Definition 4.27]{Wak5}.

\SSP
\bde[cf. ~\cite{Wak5}, Definition 4.37, (i)] \label{Def5}
An element $D^\clubsuit$ of $H^0 (X, \Sigma^{-1}(1))$ is said to be  an {\bf $(n, \vartheta)$-projective connection} on $\msX$ if the log connection $\mr{det} (D^{\clubsuit \Rightarrow \diamondsuit})$ on $\mr{det}(\mcF_\varTheta)$ induced by $D^{\clubsuit \Rightarrow \diamondsuit}$ commutes with $\nabla_\vartheta$ via \eqref{Eq22}. 
\ede
\SSP

Let $i \in \{1, \cdots, r \}$, and 
 write $\overline{\partial}$ for the section $\sigma_i^* (\mcT)$ corresponding to $1$ under  the dual of the residue isomorphism $\mr{Res} : \sigma^*_i (\Omega) \xrightarrow{\sim} k$.
 Then, $\sigma_i^* (\mcD_{})$ has a natural identification $\sigma_i^* (\mcD_{}) = k [\overline{\partial}]$ with the polynomial ring $k [\overline{\partial}]$ in $\overline{\partial}$.
Moreover, we have the following composite isomorphisms:
\begin{align} \label{Eq29}
\sigma_i^* (\mcD {\it iff}_{\vartheta, \leq n}) 
& \xrightarrow{\sim}
\mr{Hom}_k (\sigma^*_i (\varTheta^\vee), \sigma_i^* (\Omega)^{\otimes n} \otimes  \sigma_i^*(\varTheta^\vee) \otimes \sigma_i^* (\mcD_{\leq n}))
\\
& \xrightarrow{\sim}
\mr{Hom}_k (\sigma_i^* (\varTheta^\vee), \sigma_i^* (\varTheta^\vee) \otimes k [\overline{\partial}]_{\leq n}) \notag \\
& \xrightarrow{\sim}
k [\overline{\partial}]_{\leq n},
\end{align}
where $k [\overline{\partial}]_{\leq n} := \left\{ h \in k [\overline{\partial}] \, | \, \mr{deg} (h) \leq n\right\}$.
If $D^\clubsuit$ is an $(n, \vartheta)$-projective connection on $\msX$,
then, for each $i=1, \cdots, r$,
there exists a unique multiset 
\begin{align}
a_{i}(D^\clubsuit) \ \left(\text{or} \ a_{\sigma_i}(D^\clubsuit)\right) := [a_{i, 1} (D^\clubsuit), \cdots, a_{i, n} (D^\clubsuit)]
  \in \mfS_n \backslash k^{n}
\end{align}
such that the image of $\sigma_i^* (D^\clubsuit)$ via  \eqref{Eq29} coincides with 
$\prod_{j=1}^n (\overline{\partial} - a_{i, j} (D^\clubsuit))$.
Since $\mr{det}(D^{\clubsuit \Rightarrow \diamondsuit}) = \nabla_\vartheta$ via \eqref{Eq22},
the residue $\mr{Res}_{\sigma_i} (\nabla_\vartheta) \in k$ of $\nabla_\vartheta$ at  $\sigma_i$ coincides with  $(-1) \cdot \sum_{j=1}^n a_{i, j} (D^\clubsuit)$.

\SSP
\bde \label{Def65}
The multiset $a_i (D^\clubsuit)$ is referred to as the {\bf (characteristic) exponent} of $D^\clubsuit$ at $\sigma_i$.
Also, if $\rho := (\rho_1, \cdots, \rho_r) \in \mfc (k)^{r}$ denotes the $r$-tuple  determined by   $(a_1 (D^\clubsuit), \cdots a_r (D^\clubsuit))$ via the quotient $\mfS_n \backslash k^{n} \twoheadrightarrow \mfc (k)$,
then 
$\rho$ (resp., each  $\rho_i$) is referred to as the {\bf radii} of $D^\clubsuit$ (resp., the {\bf radius} of $D^\clubsuit$ at $\sigma_i$).
When $r = 0$, any $(n, \vartheta)$-projective connection is said  to  be of radii $\emptyset$, as well as of exponent $\emptyset$.
\ede
\SSP

We fix  $\rho := (\rho_1, \cdots, \rho_r) \in \mfc (k)^{r}$.
If $D^\clubsuit$ is an $(n, \vartheta)$-projective  connection of radii $\rho$,
then the characteristic polynomial of the residue matrix $\mr{Res}_{\sigma_i} (D^{\clubsuit \Rightarrow \diamondsuit})$ of $D^{\clubsuit \Rightarrow \diamondsuit}$ at $\sigma_i$ (for each $i =1, \cdots , r$)
coincides with $\prod_{j=1}^n (\overline{\partial} - a_{i, j} (D^\clubsuit))$.
That is to say, the $(\mr{GL}_n, \vartheta)$-oper  $D^{\clubsuit \Rightarrow \diamondsuit}$  turns out to be of radii $\rho$.
The resulting assignment $D^\clubsuit \mapsto D^{\clubsuit \Rightarrow \diamondsuit}$ determines a bijective correspondence 
\begin{align} \label{Eq33}
\begin{pmatrix} \text{the set of $(n, \vartheta)$-projective} \\ \text{connections on $\msX$ of radii $\rho$} \end{pmatrix} \xrightarrow{\sim} 
\begin{pmatrix} \text{the set of isomorphism classes} \\ \text{of $(\mr{GL}_n, \vartheta)$-opers on $\msX$ of radii $\rho$}\end{pmatrix}
\end{align}
(cf. ~\cite[Theorem 4.49]{Wak5}).

Next, let us take a global section $D^\clubsuit$ of $\mcD {\it iff}_{\!\vartheta, \leq n}$, and, as mentioned above, we regard it as an $n$-th order differential operator $\varTheta^\vee \rightarrow \Omega^{\otimes n} \otimes \varTheta^\vee$.
Under this interpretation, the kernel  $\mr{Ker}(D^\clubsuit)$ naturally acquires  an $\mcO_{X^{(1)}}$-module structure  via the underlying homeomorphism of $F_{X/k}$.
We say  that $D^\clubsuit$ {\bf has a full set of root functions} if $\mr{Ker} (D^\clubsuit)$ forms  a vector bundle on $X^{(1)}$ of rank $n$ (cf. ~\cite[Definition 4.64]{Wak5}).

According to 
 ~\cite[Proposition 4.65]{Wak5},
a given $(n, \theta)$-projective connection has a full set of root functions if and only if the corresponding  $(\mr{GL}_n,  \vartheta)$-oper is dormant.
Therefore, 
 the correspondence \eqref{Eq33} restricts to a bijection
\begin{align} \label{Eq37}
\begin{pmatrix} \text{the set of $(n, \vartheta)$-projective} \\ \text{connections on $\msX$ of radii $\rho$} \\ \text{with a full set of root functions}\end{pmatrix} \xrightarrow{\sim} 
\begin{pmatrix} \text{the set of isomorphism classes} \\ \text{of dormant $(\mr{GL}_n, \vartheta)$-opers } \\ \text{on $\msX$ of radii $\rho$}\end{pmatrix}.
\end{align}

\LSP
\subsection{Almost non-logarithmic extensions  of  local dormant opers} \label{SS36}

In this subsection,  we work  with  $(\mr{GL}_n, \vartheta)$-opers in a local setting, which 
 will be applied in the proof of Theorem \ref{Prop90}.

Let us write   $U := \mr{Spec} (k[\![t]\!])$, which is  equipped with 
a distinguished point  $\sigma_0$
 determined by ``$t=0$".
  We  set $U^o := \mr{Spec} (k (\!(t)\!)) \left(= U \setminus \{ \sigma_0 \}  \right)$, $\msU := (U, \{ \sigma_0 \})$, and   $\Omega :=  \left( \varprojlim_{m}\Omega_{\mr{Spec}(k[t]/(t^m))/k} \right)(\sigma_0)$.
Note that all  the definitions and constructions  discussed above remain valid when  the underlying pointed  curve $\msX$  is replaced with $\msU$.
In particular, we  may speak of  $n$-theta characteristics of $\msU$ and  $(\mr{GL}_n, \vartheta)$-opers  on $\msU$ (for an $n$-theta characteristic $\vartheta$), etc.

Let us  now fix  an $n$-theta characteristic $\vartheta := (\varTheta, \nabla_\vartheta)$ of $\msU$ such that $\nabla_\vartheta$ has vanishing $p$-curvature.
(Thus, we obtain the associated filtered vector bundle $\mcF_\varTheta := \mcD_{\leq n-1} \otimes \varTheta$, as in the global setting.)
Suppose that  there exists an integer $d$ in  $\{1, \cdots, p-n+1 \}$ satisfying  $\mr{Res}_{\sigma_0} (\nabla_\vartheta) = (-1) \cdot \sum_{j=0}^{n-2} (d+ j)$.

Let $\nabla^\diamondsuit$ be 
 a dormant $(\mr{GL}_n, \vartheta)$-oper   on $\msU$ 
 whose radius at $\sigma_0$ is given by $\rho := \LL 0, d, d +1, \cdots, d +n-2 \RR$.
Denote by $\nabla^{\diamondsuit \vee}$ the log connection on the dual bundle  $\mcF_\varTheta^\vee$  induced by $\nabla^\diamondsuit$.
The residue matrix of  $\nabla^{\diamondsuit \vee}$ at $\sigma_0$ is conjugate to the diagonal matrix with diagonal entries  $0, d', d' +1, \cdots, d' + n-2$, where $d' := p-(d+ n -2)$.
The kernel  $ \mr{Ker}(\nabla^{\diamondsuit \vee})$ naturally acquires  an $\mcO_{U^{(1)}}$-module structure  via the underlying homeomorphism of $F := F_{U/k}$.
The inclusion $\mr{Ker}(\nabla^{\diamondsuit \vee})  \hookrightarrow F_{*}(\mcF_{\varTheta}^\vee)$ corresponds, via the adjunction relation ``$F^* (-) \dashv F_* (-)$", to an $\mcO_U$-linear  morphism $\mcG  \rightarrow \mcF_\varTheta^\vee$, where $\mcG := F^* (\mr{Ker}(\nabla^{\diamondsuit \vee}))$.
This morphism is injective and becomes an isomorphism when restricted over $U^o$.
Using this injection, we  regard  $\mcG$ as an $\mcO_{U}$-submodule of $\mcF_\varTheta^\vee$.

Denote by  $\breve{\nabla}_\mcG$  the canonical  (non-logarithmic) connection on  $\mcG$,  as introduced   in  ~\cite[Theorem 5.1]{Kat1}.
This induces  a log connection $\nabla_\mcG$ on $\mcG$, which 
commutes with $\nabla^{\diamondsuit \vee}$ via the inclusion $\mcG \hookrightarrow \mcF_\varTheta^\vee$. 
For each $j=0, \cdots, n$, we define 
$\mcG^j := \mcG \cap \mcF_\varTheta^{\vee j}$,
where $\mcF_\varTheta^{\vee j} :=  (\mcF_\varTheta/\mcF_\varTheta^{n-j})^\vee \left(\subseteq \mcF_\varTheta^{\vee} \right)$.
Then, the collection $\{ \mcG^j \}_{j=0}^n$ forms an $n$-step decreasing filtration on $\mcG$ whose subquotients are line bundles.
According to (an argument similar to the proof of) ~\cite[Proposition 8.8, (i)]{Wak5},
the cokernel of the morphism
$\mcG^j/\mcG^{j+1}
\xrightarrow{}
\mcF_\varTheta^{\vee j}/ \mcF_\varTheta^{\vee j+1}$ 
induced by the inclusion $\mcG \hookrightarrow \mcF_\varTheta^\vee$ is isomorphic to the zero  sheaf
when $j=0$, and to 
  $\mcO / (t^{d' + j-1})$ when $j > 0$.

In this subsection, denote by $\breve{\mcD}$ (resp., $\breve{\Omega}$) the sheaf of {\it non-logarithmic} differential operators (resp.,  {\it non-logarithmic} $1$-forms) on $U$.
For each $j \in \mbZ_{\geq 0}$, let  $\breve{\mcD}_{\leq j}$ denote  the subsheaf of $\breve{\mcD}$ consisting of differential operators of order $\leq j$.
 For $j=1, \cdots, n-1$, the morphism $\mcF_\varTheta^{\vee j}/\mcF_\varTheta^{\vee j+1} \rightarrow \Omega \otimes (\mcF_\varTheta^{\vee j-1}/\mcF_\varTheta^{\vee j})$ induced by $\nabla^{\diamondsuit \vee}$ is an isomorphism since $\nabla^\diamondsuit$ defines a $(\mr{GL}_n, \vartheta)$-oper.
From the equality   $\nabla^{\diamondsuit \vee} |_{\mcG} = \nabla_\mcG$ together with the above argument, we deduce  that
the cokernel of the morphism 
$\mcG^j/\mcG^{j+1} 
\rightarrow \breve{\Omega} \otimes (\mcG^{j-1}/\mcG^j)$
induced by $\breve{\nabla}_\mcG$ is isomorphic to $\mcO_U/(t^{d'-1})$
when $j=1$, and isomorphic  to  $0$
  when $j > 1$.
It follows that 
the composite
\begin{align}
\kappa : 
\breve{\mcD}_{\leq n -1} \otimes   \mcG^{n-1} \xrightarrow{\mr{inclusion}} 
\breve{\mcD}_{} \otimes   \mcG \xrightarrow{\breve{\nabla}_\mcG}
\mcG,
\end{align}
restricts to an isomorphism 
$\breve{\mcD}_{\leq n -2} \otimes    \mcG^{n-1} \xrightarrow{\sim}
\mcG^1$.
Moreover, $\kappa$
 is  an isomorphism over $U^o$, and  its cokernel  is of  length $d' -1 \left(= p- d -n +1  \right)$.
Therefore,  the square 
\begin{align} \label{Eq155}
\vcenter{\xymatrix@C=46pt@R=36pt{
( \mcG/\mcG^1)^\vee
 \ar[r] \ar[d]_-{\mr{inclusion}} & ((\breve{\mcD}_{\leq n-1}/\breve{\mcD}_{\leq n-2}) \otimes \mcG^{n-1})^\vee \ar[d]^-{\mr{inclusion}} \\
\mcG^\vee  \ar[r]_-{\kappa^\vee} & (\breve{\mcD}_{\leq n-1} \otimes \mcG^{n-1})^\vee
    }}
\end{align}
forms a  pushout  diagram, where the upper horizontal arrow is the dual of the natural morphism 
\begin{align}
\left( (\breve{\mcD}_{\leq n-1} \otimes \mcG^{n-1})/ (\breve{\mcD}_{\leq n-2} \otimes \mcG^{n-1})=\right) 
(\breve{\mcD}_{\leq n-1}/\breve{\mcD}_{\leq n-2}) \otimes \mcG^{n-1} \rightarrow \mcG/\mcG^1
\end{align}
induced by  $\kappa$.

We now set
 \begin{align}
 \mcH :=  (\breve{\mcD}_{\leq n-1} \otimes \mcG^{n-1})^\vee \ \ \  \text{and}  \ \ \ \mcH^j := ((\breve{\mcD}_{\leq n-1}/\breve{\mcD}_{\leq j-1})\otimes \mcG^{n-1})^{\vee} \  (j=0, \cdots, n)
 \end{align}
  for convenience.
We regard each  $\mcH^j$   as a subbundle of $\mcH$, and  regard  $\mcG^\vee$ as an $\mcO_X$-submodule of $\mcH$ via the dual $\kappa^\vee$ of $\kappa$.
The composite injection
\begin{align} \label{Eq441}
\mcF_\varTheta \xrightarrow{} \mcG^\vee \xrightarrow{\kappa^\vee} \mcH,
\end{align}
where the first arrow denotes the dual of the inclusion $\mcG \hookrightarrow \mcF_\varTheta^\vee$,
restricts to $\mcF_\varTheta^j \hookrightarrow \mcH^j$ for each $j$, and  is an isomorphism over $U^o$.

Since the upper horizontal arrow in the above square diagram  is an inclusion between line bundles,
Note  that $\nabla_\mcG$ extends naturally  to a log connection $\nabla_\mcH$ on 
$\mcH$.
In fact, by  the upper horizontal arrow in  \eqref{Eq155},
one can identify $\mcH^{n-1}$  with $\frac{1}{t^{d'-1}} \cdot (\mcG/\mcG^1)^\vee$.
Then,
$\nabla_\mcH$  is defined in such  a way that
$\nabla_\mcH (v) := \nabla_\mcG (v)$ for  $v \in \mcG^\vee$ and 
\begin{align}
\nabla_\mcH \left(\frac{1}{t^{d'-1}} \cdot u\right) := d \left( \frac{1}{t^{d'-1}} \right) \otimes u + \frac{1}{t^{d'-1}} \cdot  \nabla_\mcG (u)
\left(=   \frac{dt}{t} \otimes  \left( \frac{1-d'}{t^{d'-1}} \cdot u\right) +  \frac{1}{t^{d'-1}}  \cdot  \nabla_\mcG (u) \right)
\end{align}
 for $u \in (\mcG/\mcG^1)^\vee$.
We thus obtain a collection of data 
\begin{align} \label{Eq159}
(\mcH, \nabla_\mcH, \{ \mcH^j \}_{j=0}^n).
\end{align}
This collection satisfies the following properties:
\begin{itemize}
\item
The log connection $\nabla^\diamondsuit$ and the filtration $\{ \mcF_\varTheta^j \}_j$ are compatible, via 
the natural inclusion $\mcF_\varTheta \hookrightarrow \mcH$, with $\nabla_\mcH$ and $\{ \mcH^j \}_{j}$, respectively.
Moreover,   we have the identity 
\begin{align}
(\mcH, \nabla_\mcH, \{ \mcH^j \}_{j=0}^n) |_{U^o} = (\mcF_\varTheta, \nabla^\diamondsuit, \{ \mcF_\varTheta^j \}_{j=0}^n) |_{U^o};
\end{align}
\item
For each $j=0, \cdots, n-1$, 
the cokernel of the injection  $\mcF_\varTheta^j/\mcF_\varTheta^{j+1} \hookrightarrow \mcH^j /\mcH^{j+1}$ is of  length $p-d - j$;
\item
The exponent of $\nabla_\mcH$ coincides with $[0, 0, \cdots, 0, d+ n-1]$.
\end{itemize}

\SSP
\bde \label{Def50}
We shall refer to the collection $(\mcH, \nabla_\mcH, \{ \mcH^j \}_{j})$ as the {\bf almost non-logarithmic extension} of $(\mcF_\varTheta, \nabla^\diamondsuit, \{ \mcF_\varTheta^j \}_{j})$.
\ede
\SSP

Next, let $\nabla_\circ^\diamondsuit$ and $\nabla_\bullet^\diamondsuit$ be dormant $(\mr{GL}_n, \vartheta)$-opers on $\msU$ of radii $\rho$.
For each $\Box \in \{ \circ, \bullet \}$, we denote by $(\mcH_\Box, \nabla_{\mcH, \Box}, \{ \mcH_\Box^j \}_{j=1}^n)$ the almost non-logarithmic  extension  of $(\mcF_\varTheta, \nabla_\Box^\diamondsuit, \{ \mcF_\varTheta^j \}_{j=1}^n)$.
Define
\begin{align}
\mr{Isom} (\nabla_\circ^\diamondsuit, \nabla_\bullet^\diamondsuit) \ \left(\text{resp.,} \  \mr{Isom} (\nabla_{\mcH, \circ}, \nabla_{\mcH, \bullet}) \right)
\end{align}
to be  the set of filtered automorphisms 
$\mcF_\varTheta \xrightarrow{\sim} \mcF_\varTheta$ (resp., filtered isomorphisms $\mcH_\circ \xrightarrow{\sim} \mcH_\bullet$) via which $\nabla_\circ^\diamondsuit$  (resp., $\nabla_{\mcH, \circ}$) commutes with $\nabla_\bullet^\diamondsuit$ (resp., $\nabla_{\mcH, \bullet}$).  
Given  an element $h$ of $\mr{Isom} (\nabla_\circ^\diamondsuit, \nabla_\bullet^\diamondsuit)$,
 it follows from the functorial construction of  almost non-logarithmic extension  that
$h$ uniquely extends to an element $h_\mcH$ of $\mr{Isom} (\nabla_{\mcH, \circ}, \nabla_{\mcH, \bullet})$.

\SSP
\bpr \label{PRop532}
The map of sets
\begin{align} \label{Eq901}
\mr{Isom} (\nabla_\circ^\diamondsuit, \nabla_\bullet^\diamondsuit) \xrightarrow{} \mr{Isom} (\nabla_{\mcH, \circ}, \nabla_{\mcH, \bullet})
\end{align}
given by 
the resulting assignment $h \mapsto h_\mcH$ is bijective.
\epr
\begin{proof}
By the composite injection  \eqref{Eq441},
  $\mcF_\varTheta$ can be regarded as a subsheaf of $\mcH_\circ$, as well as of  $\mcH_\bullet$.
The injectivity of  \eqref{Eq901} follows immediately  from this observation.

To prove the surjectivity, let us take an element $\widetilde{h}$ of  $\mr{Isom} (\nabla_{\mcH, \circ}, \nabla_{\mcH, \bullet})$.
For each $\Box \in \{ \circ, \bullet \}$,
denote by $\mcG_\Box$ (resp., $\nabla_{\mcG, \Box}$; resp., $\kappa_\Box$) the sheaf  ``$\mcG$" (resp., the log connection ``$\nabla_\mcG$"; resp.,  the morphism ``$\kappa$") obtained by applying the construction described above  to $\nabla_\Box^\diamondsuit$.
The inclusion $\mr{Ker} (\nabla_{\mcH, \Box}) \hookrightarrow \mcH_\Box$ induces an $\mcO_U$-linear morphism $\iota_\Box : F^* (\mr{Ker} (\nabla_{\mcH, \Box})) \rightarrow \mcH_\Box$.
This morphism is injective,  and the dual of the pair $(F^* (\mr{Ker} (\nabla_{\mcH, \Box})), \iota_\Box)$ coincides  with 
$(\mcG_\Box, \kappa_\Box)$.
Moreover,
 $\nabla_{\mcG,  \Box}$ corresponds to  
the canonical  connection resulting from ~\cite[Theorem 5.1]{Kat1}
 under the natural identification $F^* (\mr{Ker}(\nabla_{\mcH, \Box}))^\vee = F^* (\mr{Ker}(\nabla_{\mcH, \Box})^\vee)$.
This  connection commutes with $\nabla^{\diamondsuit \vee}_\Box$ via
  the inclusion $\mcG_\Box \hookrightarrow \mcF_\varTheta^\vee$.
On the other hand, 
since 
$\widetilde{h}$ preserves the log connection, it yields a filtered  isomorphism $h' : (\mcG_\bullet, \nabla_{\mcG, \bullet}) \xrightarrow{\sim} (\mcG_\circ, \nabla_{\mcG, \circ})$.
Let  $\mcG_{\Box, +}$ denote the pushout of the inclusions  $\mcG^{n-1}_\Box \hookrightarrow \mcG_\Box$ and $\mcG^{n-1}_\Box \hookrightarrow \mcF_\varTheta^{\vee n-1}$.
Then, $h'$ extends to  an isomorphism of vector bundles $h'_+ : \mcG_{\bullet, +} \xrightarrow{\sim} \mcG_{\circ, +}$.
Since $\mcF^\vee_{\varTheta}$ is generated by sections of $\mcF^{\vee n-1}_\varTheta$ as  a logarithmic  flat bundle,
the isomorphism  $h'_+$ further  extends to an isomorphism  $(\mcF_\varTheta^\vee, \nabla_\bullet^{\diamondsuit \vee}, \{ \mcF_\varTheta^{\vee j} \}_j) \xrightarrow{\sim} (\mcF_\varTheta^\vee, \nabla_\circ^{\diamondsuit \vee}, \{ \mcF_\varTheta^{\vee j} \}_j)$.
Taking the dual of this isomorphism defines an element $h$ of $\mr{Isom}(\nabla_\circ^\diamondsuit, \nabla_\bullet^\diamondsuit)$, and one can verify the identity  $h_\mcH = \widetilde{h}$.
This completes the proof of the desired surjectivity  of \eqref{Eq901}. 
\end{proof}

\LSP
\subsection{$n$-theta characteristics of a $3$-pointed projective line} \label{SS12}

The remainder of this section  focuses  on the case where $\msX$ is taken to be the  $3$-pointed projective line 
$\msP := (\mbP, \{ [0], [1], [\infty]\})$, where $\mbP := \mr{Proj}(k[s, t])$ and 
  for each $\lambda \in k \sqcup \{ \infty \}$ we denote by   $[\lambda]$  the corresponding  $k$-rational point of $\mbP$.
  The  set $\{[0], [1], [\infty] \}$ is 
   considered as  the ordered  set of $3$ marked points.

We define the local coordinates   $x := \frac{s}{t}$, $y :=  \frac{1}{x}\left(= \frac{t}{s} \right)$, and $z := x-1$.
For $w \in \{ x, y, z \}$, we write  
$U_w$ for  the formal neighborhood of the point $w = 0$ in $\mbP$, and 
 set  $U_w^o := U_w \setminus \{ w = 0\}$.
We also define the derivation  $\delta_w := w  \frac{d}{dw}$, which acts locally on  $\mcO_\mbP$, and restricts to a derivation on $\mcO_{U_w}$.
The assignment  $v \mapsto v \cdot \delta_w$ gives   an identification $\eta_w : \mcO_{U_w} \xrightarrow{\sim} \mcT |_{U_w}$.

Let us now fix 
two collections $\alpha := (\alpha_1, \cdots, \alpha_n) \in k^{n}$, $\beta := (\beta_1, \cdots, \beta_{n-1}) \in k^{n-1}$.
The pair $(\alpha, \beta)$ determines a log connection $\nabla_{\alpha, \beta}$ on  $\mcT^{\otimes \frac{n(n-1)}{2}} |_{U_x}$ 
 expressed as 
\begin{align}
\nabla_{\alpha, \beta} = d +  \frac{dx}{x} \otimes \frac{x \cdot \sum_{j=1}^n \alpha_j - \sum_{j=1}^{n-1}(\beta_j -1)}{x-1}
\end{align}  
under the identification 
$\mcO_{U_x} = \mcT^{\otimes \frac{n(n-1)}{2}} |_{U_x}$ 
given by
$\eta^{\otimes \frac{n(n-1)}{2}}_x$.
This expression  extends  to a log connection on $ \mcT^{\otimes \frac{n(n-1)}{2}} |_{\mbP \setminus \{[1], [\infty] \}}$,
 and its restriction to $U_y^o$ (resp., $U_z^o$) is given by  
\begin{align}
\nabla_{\alpha, \beta} |_{U_y^o} = d + \frac{dy}{y} \otimes \frac{\sum_{j=1}^n \alpha_j - y \cdot \sum_{j=1}^{n-1} (\beta_j -1)}{y-1} \hspace{28mm}  \\
\left(\text{resp.,} \ \nabla_{\alpha, \beta} |_{U_z^o} = d + \frac{dz}{z} \otimes \frac{(z+1) \sum_{j=1}^n \alpha_j  - \sum_{j=1}^{n-1}(\beta_j -1) +\frac{n(n-1)}{2}}{z+1} \right)
\end{align}
under the identification $\mcO_{U_y^o} = \mcT^{\otimes \frac{n (n-1)}{2}} |_{U_y^o}$ (resp., $\mcO_{U_z^o} = \mcT^{\otimes \frac{n(n-1)}{2}} |_{U_z^o}$) given by $\eta_y^{\otimes \frac{n(n-1)}{2}}$ (resp., $\eta_z^{\otimes \frac{n(n-1)}{2}}$).
Hence, $\nabla_{\alpha, \beta}$ defines a global log connection on $\mcT$, and the pair 
\begin{align} \label{Eq44d1}
\vartheta_{\alpha, \beta} := (\mcO_{\mbP}, \nabla_{\alpha, \beta})
\end{align}
 specifies an $(n, \theta)$-theta characteristic of $\msP$.

\SSP
\bpr \label{Prop22}
Let us keep the above notation.
Then,  both  $\sum_{j=1}^{n} \alpha$ and $\sum_{j=1}^{n-1} \beta_j$ belong to $\mbF_p$ if and only if $\nabla_{\alpha, \beta}$ has vanishing $p$-curvature.
\epr
\begin{proof}
Note that $\nabla_{\alpha, \beta} = d + Q$, where 
\begin{align}
Q := \left(-\sum_{j=1}^n \alpha + \sum_{j=1}^{n-1} (\beta_j -1)\right) \cdot \frac{d (x-1)}{x-1} - \left(\sum_{j=1}^{n-1}(\beta_j -1)\right) \cdot \frac{dx}{x}.
\end{align}
According to  ~\cite[Corollary 7.1.3]{Kat2},
the log connection $\nabla_{\alpha, \beta}$ has vanishing $p$-curvature if and only if $Q$ is invariant under the Cartier operator on $\Omega$, in the sense of ~\cite[(7.1.3.2)]{Kat2}.
Since both $\frac{d (x-1)}{x-1}$ and $\frac{dx}{x}$ are invariant under the Cartier  operator, 
the latter condition of  the desired   equivalence translates into the requirement  that both $\sum_{j=1}^{n} \alpha_j$ and $\sum_{j=1}^{n-1} \beta_j$ are invariant under  the Frobenius endomorphism of $\mr{Spec} (k)$.
This holds if and only if these sums  lie in $\mbF_p$.
This completes the proof of the assertion.
\end{proof}

\LSP
\subsection{Projective connections arising from  generalized hypergeometric operators} \label{SS32}

Let $(\alpha, \beta)$ be as above, and define  the multisets $a_1^{\alpha, \beta}$, $a_2^{\alpha, \beta}$, and  $\alpha_3^{\alpha, \beta}$ as follows:
\begin{align} \label{Eq333}
a^{\alpha, \beta}_{1} &:= [0, 1- \beta_1, 1- \beta_2, \cdots, 1- \beta_{n-1}], \\
a_{2}^{\alpha, \beta}  &:= [0, 1, 2, \cdots, n-2, \sum_{j=1}^{n-1} \beta_i - \sum_{j=1}^n \alpha_i], \notag \\
a_{3}^{\alpha, \beta} &:= [\alpha_1, \alpha_2, \alpha_3, \cdots, \alpha_n].  
\end{align}
Denote by $\rho^{\alpha, \beta}_1$, $\rho_2^{\alpha, \beta}$, and $\rho_3^{\alpha, \beta}$ the elements of $\mfc (k)$ determined by $a_1^{\alpha, \beta}$, $a_2^{\alpha, \beta}$, and $a^{\alpha, \beta}_3$, respectively.
We also set $a^{\alpha, \beta} := (a_1^{\alpha, \beta}, a_2^{\alpha, \beta}, a_3^{\alpha, \beta})$ and $\rho^{\alpha, \beta} := (\rho_1^{\alpha, \beta}, \rho_2^{\alpha, \beta}, \rho_3^{\alpha, \beta})$.

Now, consider the operator  $D_{\alpha, \beta}^\clubsuit  := \frac{dx^{\otimes n}}{x^n (1-x)} \otimes D_{\alpha, \beta}$ (cf. \eqref{Eq6} for the definition of $D_{\alpha, \beta}$), which defines 
  an $n$-order linear differential operator 
from $\mcO_{\mbP \setminus \{ [1], [\infty]  \}}$ to $\Omega^{\otimes n} |_{\mbP \setminus \{ [1], [\infty]\}}$.
Explicitly, we write 
\begin{align}
D^\clubsuit_{\alpha, \beta}  = \frac{dx^{\otimes n}}{x^n (1-x)} \otimes \left( \delta_x \cdot \prod_{j=1}^{n-1} (\delta_x + \beta_j -1) - x \cdot  \prod_{j=1}^n (\delta_x + \alpha_j) \right).
\end{align}
As discussed in ~\cite[Remark 4.12]{Wak5},
this operator  may be regarded as a  global section of $\mcD {\it iff}_{\!\vartheta_{\alpha, \beta}, \leq n}$ over $\mbP \setminus \{[1], [\infty] \}$, and it is mapped to $1$ via the symbol map $\Sigma$.

The following assertion is well-known at least  in the complex analytic setting (cf., e.g.,  ~\cite[Section 2]{BeHe}).
However, we are not aware of any reference that treats this fact algebraically, particularly including the  computation of the exponent at $x=1$ {\it in  positive characteristic}.
Thus, we include an algebraic  proof here for the reader's benefit.

\SSP
\bpr \label{Lem3}
Let us retain the notation established above.
Then,
 $D^\clubsuit_{\alpha, \beta}$ extends to an $n$-order linear differential operator $\mcO_{\mbP} \rightarrow \Omega^{\otimes n}$, which defines    
an $(n, \theta_{\alpha, \beta})$-projective connection on $\msP$.
Moreover, 
the  exponents of $D^\clubsuit_{\alpha, \beta}$ at $[0]$, $[1]$, and $[\infty]$ are, respectively,  given by 
\begin{align}
a_{[0]} (D_{\alpha, \beta}^\clubsuit) =a_1^{\alpha, \beta},
\hspace{10mm}
a_{[1]} (D_{\alpha, \beta}^\clubsuit) =a_2^{\alpha, \beta},
\hspace{10mm}
a_{[\infty]} (D_{\alpha, \beta}^\clubsuit) =a_3^{\alpha, \beta}.
\end{align}
\epr
\begin{proof}
On the formal neighborhood $U_x$,   we have
\begin{align}
D^\clubsuit_{\alpha, \beta} |_{U_x}  = \left(\frac{dx}{x}\right)^{\otimes n} \otimes \left(\frac{1}{1-x} \cdot \delta_x \cdot \prod_{j=1}^n (\delta_x + \beta_j -1)  - \frac{x}{1-x} \cdot  \prod_{j=1}^n (\delta_x + \alpha_j) \right). 
\end{align}
 From this expression,  it follows   that  the exponent at $x= 0$ (i.e., at the point $[0]$) coincides with 
 $a_1^{\alpha, \beta}$.
Next, the restriction of $D^\clubsuit_{\alpha, \beta}$ 
to  $U_y^o$ can be described as 
\begin{align}
D^\clubsuit_{\alpha, \beta} |_{U_y^o}  = \left(\frac{dy}{y} \right)^{\otimes n} \otimes \left( \frac{-y}{1-y} \cdot \delta_y  \cdot \prod_{j=1}^{n-1} (\delta_y -\beta_j +1) + \frac{1}{1-y} \cdot \prod_{j=1}^n (\delta_y - \alpha_j)\right).
\end{align}
This expression shows that $D^\clubsuit_{\alpha, \beta} |_{U_y^o}$ extends to an $n$-th order  linear differential operator $\mcO_{U_y} \rightarrow \Omega^{\otimes n} |_{U_y}$,
 and that the exponent at $y =0$ (corresponding to $x = \infty$, i.e.,  the point $[\infty]$) coincides with $a_3^{\alpha, \beta}$.

Finally, let us consider the restriction of $D_{\alpha, \beta}^\clubsuit$ to $U_z$.
For convenience, we set $\beta_n :=1$.
Observe that,
for any $\nu \in k$ and $j \in \mbZ_{\geq 0}$,
the  following identity holds: 
\begin{align}
z^j \cdot  \left(\delta_z + \frac{1}{z} \cdot \delta_z + \nu \right) = \left(\delta_z + \frac{1}{z} \cdot \delta_z + \nu - j \left(1 + \frac{1}{z} \right) \right) \cdot z^j.
\end{align}
Using this identity, we compute 
the following sequence of equalities:
\begin{align} \label{Eq144}
 & \ \ \  \ D_{a, b}^{\clubsuit}  |_{U_z^o} \\
& = \frac{dz^{\otimes n}}{(z+1)^nz^{n+1}} \otimes z^n \cdot \left(-\prod_{j=1}^{n} \left(\delta_z + \frac{1}{z} \cdot \delta_z + \beta_j -1\right) + (z +1) \cdot \prod_{j=1}^n \left(\delta_z + \frac{1}{z} \cdot \delta_z + \alpha_j\right) \right)  \notag \\
& =  \frac{dz^{\otimes z}}{(z+1)^nz^{n+1}} \otimes \Biggl( - \prod_{j=1}^n z \cdot  \left(\delta_z + \frac{1}{z} \cdot \delta_z + \beta_j -1 -(j-1) \left(1+ \frac{1}{z} \right) \right) \notag \\
& \hspace{10mm}  + (z+1) \cdot  \prod_{j=1}^n z \cdot \left(\delta_z + \frac{1}{z} \cdot \delta_z + \alpha_j -(j-1) \left(1 + \frac{1}{z} \right) \right) \Biggl)
\notag \\
& = \left(\frac{dz}{z} \right)^{\otimes n} \otimes \frac{1}{(z+1)^nz}
\cdot \Biggl(-\prod_{j=1}^n \left(\delta_z -(j-1) + z (\delta_z + \beta_j -j) \right) \notag \\
& \hspace{10mm} + (z+1) \cdot \prod_{j=1}^n \left(\delta_z -(j-1) + z (\delta_z + \alpha_j-j +1)\right)\Biggl), \notag 
\end{align}
where, for non-commuting  differential operators $Q_1, \cdots, Q_n$,  the notation $\prod_{j=1}^n Q_j$ denotes their ordered composition $Q_n \circ Q_{n-1} \circ  \cdots \circ  Q_1$.
Focusing  on the final expression above,
we note that the leading terms (with respect to $\delta_z$) of the linear differential operators
\begin{align}
\prod_{j=1}^n \left(\delta_z -(j-1) + z (\delta_z + \beta_j -j) \right), \hspace{10mm}
(z+1) \cdot \prod_{j=1}^n \left(\delta_z -(j-1) + z (\delta_z + \alpha_j-j +1)\right)
\end{align}
 coincide  modulo $z$.
Hence, $D_{a, b}^{\clubsuit}  |_{U_z^o}$  extends  to an $n$-th order  differential operator $\mcO_{U_z} \rightarrow \Omega^{\otimes n}|_{U_z}$.

Moreover, modulo $z$, the rightmost of the above sequence
can be  computed  as follows:
\begin{align}
&  \ \ \ \ (\text{RHS of \eqref{Eq144}}) \\
& =  \left(\frac{dz}{z} \right)^{\otimes n} \otimes \frac{1}{(z+1)^nz}
\cdot \Biggl( -\sum_{j=1}^n  \left(\prod_{s=j+1}^n \left(\delta_z - s+1 \right)\right) \cdot z (\delta_z +\beta_j -j) \cdot \left(\prod_{s=1}^{j-1} (\delta_z - s+1) \right) \notag \\
& \hspace{10mm}+  \sum_{j=1}^n \left( \prod_{s=j+1}^n \left(\delta_z - s+1 \right)\right) \cdot z (\delta_z + \alpha_j -j+1) \cdot \left(\prod_{s=1}^{j-1} (\delta_z - s+1) \right) \Bigg) \notag \\
& =  \left(\frac{dz}{z} \right)^{\otimes n} \otimes \frac{1}{(z+1)^nz}\cdot
 \Bigg(\sum_{j=1}^n \left( \prod_{s=j+1}^n \left(\delta_z - s+1 \right)\right) \cdot z (\alpha_j - \beta_j +1) \cdot \left(\prod_{s=1}^{j-1} (\delta_z - s+1) \right) 
 \notag \\
 & \hspace{10mm} +z \cdot \prod_{j=1}^n (\delta_z -j+1)\Bigg) \notag
  \\
 & =  \left(\frac{dz}{z} \right)^{\otimes n} \otimes \frac{1}{(z+1)^nz}\cdot \Biggl( \sum_{j=1}^n \left( z \cdot (\alpha_j -\beta_j +1) \cdot \prod_{s=j+1}^n \left(\delta_z - s+2 \right)\right) \cdot \left(\prod_{s=1}^{j-1} (\delta_z - s+1) \right) \notag
  \\
 & \hspace{10mm} + z \cdot \prod_{j=1}^n (\delta_z -j+1)\Bigg) \notag  \\
 & =   \left(\frac{dz}{z} \right)^{\otimes n} \otimes \frac{1}{(z+1)^n}\cdot 
 \left( \left(\sum_{j=1}^n (\alpha_j -\beta_j +1)\right) \cdot  \prod_{s =0}^{n-2} (\delta_z -s) + \prod_{s = 0}^{n-1} (\delta_z - s) \right)  \notag
  \\
 & = \left(\frac{dz}{z} \right)^{\otimes n} \otimes \frac{1}{(z+1)^n} \cdot \left(\delta_z - \left(\sum_{j=1}^{n-1} \beta_j - \sum_{j=1}^n\alpha_j\right)\right) \cdot \prod_{s=0}^{n-2} (\delta_z - s). \notag
 \end{align}
Therefore, the exponent at $z =0$ (corresponding to $x =1$, i.e., the point $[1]$) is given by $a_2^{\alpha, \beta}$.
Since the definition of $D_{\alpha, \beta}^\clubsuit$ implies  $\mr{det}(\nabla_{D_{\alpha, \beta}^\clubsuit}) = \nabla_{\alpha, \beta}$ via \eqref{Eq22}, we  conclude that $D^\clubsuit_{\alpha, \beta}$ defines an $(n, \vartheta_{\alpha, \beta})$-projective connection.
This  completes  the proof of this assertion. 
\end{proof}
\SSP

For each pair $(\alpha, \beta)$ as above,
the $(\mr{GL}_n, \vartheta_{\alpha, \beta})$-oper (resp., the $\mr{PGL}_n$-oper) associated to $D^\clubsuit_{\alpha, \beta}$ via  \eqref{Eq33} (resp., \eqref{Eq14} and \eqref{Eq33}) is denoted by
\begin{align}
\nabla_{\alpha, \beta}^\diamondsuit \ \left(\text{resp.,} \ \msE^\spadesuit_{\alpha, \beta} \right). 
\end{align}
In particular, Proposition \ref{Lem3} implies that $\nabla_{\alpha, \beta}^\diamondsuit$ (resp., $\msE^\spadesuit_{\alpha, \beta}$) has exponents $a^{\alpha, \beta}$ (resp., radii $\rho^{\alpha, \beta}$).

\SSP
\bde \label{Def78}
A $\mr{PGL}_n$-oper on $\msP$ is said to be {\bf of hypergeometric type}
if its radius  at  one of the marked points in $\msP$ can be represented by a multiset  of the form $[0, 1, \cdots, n-2, d]$ for some $d \in \{ n-1, n, \cdots, p-1 \}$.
\ede
\SSP

Note that a $\mr{PGL}_n$-oper is of hypergeometric type if and only if it 
  is isomorphic to 
the pull-back of 
$\msE^\spadesuit_{\alpha, \beta}$ for some pair  $(\alpha, \beta)$  via a $k$-automorphism $\eta$ of $\mbP$ satisfying $\eta (\{ [0], [1], [\infty] \}) \subseteq \{ [0], [1], [\infty] \}$.
The following assertion provides a criterion for determining  when such a  $\mr{PGL}_n$-oper is dormant.

\SSP
\bpr \label{Prop4}
Let $(\alpha, \beta)$ be as above.
Then, the following three conditions are equivalent:
\begin{itemize}
\item[(a)]
The $(\mr{GL}_n, \vartheta_{\alpha, \beta})$-oper $\nabla_{\alpha, \beta}^\diamondsuit$ is dormant;
\item[(b)]
The $\mr{PGL}_n$-oper $\msE^\spadesuit_{\alpha, \beta}$ is dormant;
\item[(c)]
The tuples  $\alpha$ and $\beta$ lie in  $\mbF_p^{n}$ and $\mbF_p^{n-1}$, respectively, and 
after possibly reordering the indices of $\{ \alpha_j \}_{j=1}^n$ and $\{ \beta_j \}_{j=1}^{n-1}$, the following inequalities hold:
\begin{align} \label{Eq182}
\widetilde{\alpha}_1 \geq \widetilde{\beta}_1 > \widetilde{\alpha}_2 \geq \widetilde{\beta}_2 > \cdots \geq \widetilde{\beta}_{n-1} > \widetilde{\alpha}_n.
\end{align}
\end{itemize}
\epr
\begin{proof}
The implication (a) $\Rightarrow$ (b) is immediate. 
Since $\mr{Ker}(D^\clubsuit_{\alpha, \beta})$ is torsion-free (or equivalently, locally free) as an  $\mcO_{\mbP^{(1)}}$-module,
the kernel   $\mr{Ker}(D_{a, b})$  has rank $n$ if and only if $D_{a, b}^\clubsuit$ has a full set of root functions.
Thus,  the equivalent (a) $\Leftrightarrow$ (c) 
 follows directly  from  Corollary \ref{Cor11}.

To prove the remaining portion,
 we suppose that
 the condition (b) holds, i.e., 
  $\msE_{\alpha, \beta}^\spadesuit$ is dormant.
 Let us take an $n$-theta characteristic  $\vartheta := (\varTheta, \nabla_\vartheta)$ of $\msP$ such that $\nabla_\vartheta$ has vanishing $p$-curvature.
 Then, $\msE^\spadesuit_{\alpha, \beta}$ corresponds to a dormant $(\mr{GL}_n, \vartheta)$-oper $\nabla^\diamondsuit$ via \eqref{Eq14}.
It follows from ~\cite[Corollary 2.10]{Oss3} (or ~\cite[Proposition-Definition 4.8]{Wak20}) that  the residue matrices  of $\nabla^\diamondsuit$ at the  points $[0]$, $[1]$, and $[\infty]$
are 
  conjugate to diagonal matrices whose entries lie in $\mbF_p$.
  Since $\nabla^\diamondsuit$ and $\nabla_{\alpha, \beta}^\diamondsuit$ differ by tensoring with a flat line bundle,
 the  multisets given by these diagonal entries for  the respective marked points coincide with  $a_1^{\alpha, \beta}$,  $a_2^{\alpha, \beta}$, and  $a_3^{\alpha, \beta}$  up to translation by scalars in $k$. 
 This implies that the sums $\sum_{j=1}^n \alpha_j$ and $\sum_{j=1}^{n-1} \beta_j$ lie in $\mbF_p$, and  $\nabla_{\alpha, \beta}$ has vanishing $p$-curvature by  Proposition \ref{Prop22}.
By the bijectivity of  \eqref{Eq14}, the $(\mr{GL}_n, \vartheta_{\alpha, \beta})$-oper  $\nabla^\diamondsuit_{\alpha, \beta}$ turns out to be dormant, i.e., the condition (a) holds.
This completes the proof of this proposition.
\end{proof}
\SSP

We denote by 
\begin{align} \label{Eq149}
\mr{Hyp}_{p, n}
\end{align}
the subset of $\mfc (\mbF_p)^{3}$ consisting of all triples $(\rho_1, \rho_2, \rho_3)$ 
for which there exists a permutation $\sigma \in \mfS_3$ and 
 a pair $(\alpha, \beta) \in \mbF_p^{n} \times \mbF_p^{n-1}$ satisfying the condition \eqref{Eq182} such that  $(\rho_{\sigma (1)}, \rho_{\sigma (2)}, \rho_{\sigma (3)}) = (\rho_{1}^{\alpha, \beta}, \rho_2^{\alpha, \beta}, \rho_3^{\alpha, \beta})$.
Proposition \ref{Prop4} asserts  that a triple $\rho \in \mfc (\mbF_p)^{3}$ belongs to $\mr{Hyp}_{p, n}$ if and only if
there exists a dormant $\mr{PGL}_n$-oper on $\msP$ of hypergeometric type whose radii are given by  $\rho$.

The following proposition establishes   a certain rigidity property  of dormant $\mr{PGL}_n$-opers of   hypergeometric type.

\SSP
\bt \label{Prop90}
Let  
$\alpha := (\alpha_1, \cdots, \alpha_n) \in \mbF_p^{n}$, $\beta := (\beta_1, \cdots, \beta_{n-1}) \in \mbF_p^{n-1}$, and suppose that 
$\rho^{\alpha, \beta} \in \mr{Hyp}_{p, n}$.
Then, any two dormant $\mr{PGL}_n$-opers on $\msP$ of radii $\rho^{\alpha, \beta}$ are isomorphic to each other.
In particular, $\msE_{\alpha, \beta}^\spadesuit$ is the unique dormant $\mr{PGL}_n$-oper on $\msP$ of radii $\rho^{\alpha, \beta}$ up to isomorphism.
\et
\begin{proof}
Set $\gamma := \sum_{j=1}^{n-1} \beta_j - \sum_{j=1}^n \alpha_n$
(hence the multiset $[0, 1, \cdots, n-2, \gamma]$  coincides with  $[0, p-\gamma, \cdots, p-\gamma + n-2]$ in $\mfc (\mbF_p)$).
 One can find an  $n$-theta characteristic $\vartheta := (\varTheta, \nabla_\vartheta)$ of $\msP$ satisfying the following conditions (cf. ~\cite[Section 4.6.4]{Wak5}):
\begin{itemize}
\item
The log connection $\nabla_\vartheta$ has vanishing $p$-curvature;
\item
The residue  $a_{[1]}$ of $\nabla_\vartheta$ at $[1]$ satisfies $a_{[1]} = \sum_{j=1}^{n-2} (p-\gamma + j)$;
\item
The underlying line bundle $\varTheta$ coincides with $\mcO_\mbP (- (p-d-n +1)[1])$.
\end{itemize}
Denote by $a_{[0]}$ (resp., $a_{[\infty]}$) the residue  of $\nabla_\vartheta$ at $[0]$ (resp., $[\infty]$).

Now, let us consider  two dormant $\mr{PGL}_n$-opers $\msE^\spadesuit_\circ$, $\msE^\spadesuit_\bullet$ on $\msP$ of radii $\rho^{\alpha, \beta}$.
For each 
$\Box \in \{ \circ, \bullet \}$,
denote by 
 $\nabla^\diamondsuit_\Box$  
the dormant $(\mr{GL}_n, \vartheta)$-oper   corresponding to $\msE^\spadesuit_\Box$ via the equivalence (a) $\Leftrightarrow$ (b) in Proposition \ref{Prop4}.
Without loss of generality,  we may assume  that $\nabla_\Box^\diamondsuit$ is normal in the sense of ~\cite[Definition 4.53]{Wak5}, i.e., $\nabla_\Box^\diamondsuit = D^{\clubsuit \Rightarrow \diamondsuit}_\Box$ for some $(n, \vartheta)$-projective connection $D^\clubsuit_\Box$ having a full set of root functions.
By the second condition listed above, the exponent  of $D_\Box^{\clubsuit}$ at $[1]$ coincides with $[0, p-\gamma, \cdots, p-\gamma + n-2]$.
Let us  define  
\begin{align}
a'_{[0]} := - \frac{a_{[0]} + \sum_{j=1}^{n-1} (1 -\beta_j)}{n} \  \left(\text{resp.,} \  a'_{[\infty]} := - \frac{a_{[\infty]} + \sum_{j=1}^n \alpha_n}{n}\right).
\end{align}
Then,  the exponent of $D_\Box^\clubsuit$ at $[0]$ (resp., $[\infty]$) is given by
$[a'_{[0]}, a'_{[0]} + 1-\beta_1, \cdots, a'_{[0]} + 1 - \beta_{n-1}]$ (resp., $[a'_{[\infty]} + \alpha_1, a'_{[\infty]} + \alpha_2, \cdots, a'_{[\infty]} + \alpha_n]$).

In this proof, let $\breve{\Omega}$ denote the sheaf of logarithmic $1$-forms on $\mbP/k$ with poles along the divisor $[0] +[\infty]$.
Write $\breve{\mcT} := \breve{\Omega}^\vee$, and write 
 $\breve{\mcD}_{\leq j}$ (for $j \in \mbZ_{\geq 0}$)  for the associated sheaf of logarithmic differential operator of order $\leq j$.
Also, we set $\breve{\varTheta} := \varTheta ((p-d -n +1) [1]) \left(\cong \mcO_{\mbP} \right)$,
 $\breve{\mcF}_\varTheta := \breve{\mcD}_{\leq n-1} \otimes \breve{\varTheta}$, and
 $\breve{\mcF}_\varTheta^j := \breve{\mcD}_{\leq n-j-1}\otimes \breve{\varTheta}$ ($j=0, \cdots, n$).
  The inclusions $\mcD \hookrightarrow \breve{\mcD}$ and $\varTheta \hookrightarrow \breve{\varTheta}$ induce  an injection 
  $\mcF_\varTheta \hookrightarrow \breve{\mcF}_\varTheta$, with respect to which  
the  $n$-step decreasing filtration  $\{ \breve{\mcF}_\varTheta^j \}_{j=0}^n$ commutes  with $\{ \mcF^j_\varTheta \}_{j=0}^n$.

The restriction of  $\nabla^\diamondsuit_\Box$ to $U_z$ (cf. \S\,\ref{SS12}) forms  a dormant $(\mr{GL}_n, \vartheta |_{U_z})$-oper on the pointed formal disc $(U_z,  \{ [1]\})$.
Using  the identification $U = U_z$ determined  by $t = z$,
we obtain  the  almost  non-logarithmic extension $(\mcH, \nabla_\mcH, \{ \mcH^j \}_{j})$ of $(\mcF_\varTheta, \nabla_\Box^\diamondsuit, \{ \mcF^j_{\varTheta} \}_j) |_{U_z}$ (cf. Definition \ref{Def50}).
The two collections   $(\mcH, \nabla_\mcH, \{ \mcH^j \}_{j})$ and  $(\mcF_\varTheta, \nabla_\Box^\diamondsuit, \{ \mcF_\varTheta^j \}_j)|_{\mbP \setminus \{ [1]\}}$ can be glued together to form  a filtered logarithmic  flat bundle on $\msP$.
By construction,
  the resulting object  takes  the form $(\breve{\mcF}_\varTheta, \breve{\nabla}_\Box^\diamondsuit, \{ \breve{\mcF}_\varTheta^j \}_{j= 0}^n)$ for some log connection $\breve{\nabla}^\diamondsuit_\Box$ on  $\breve{\mcF}_\varTheta$.

Now, consider the global section $x \frac{d}{dx} \left(= - y \frac{d}{dy} \right)$ of $\breve{\mcT}$ and the natural identification $\breve{\varTheta}  = \mcO_\mbP$.
These  together  determine a trivialization  $\breve{\mcF}_\varTheta = \mcO_\mbP^{\oplus n}$.
Under this  trivialization, 
$\breve{\nabla}_\Box^\diamondsuit$ can be expressed as $d + A$ for some $A \in H^0 (\mbP,  \Omega \otimes_k M_n (k))$, where $M_n (k)$ denotes the space of $n \times n$ matrices over $k$, regarded as $k$-linear endomorphisms of $k^{n}$. 
Since $\nabla_\Box^\diamondsuit$ is normal,
the matrix  $A$  decomposes  as $A = A_{[0]} \otimes \frac{dw}{w+1} + A_{[\infty]} \otimes \frac{dw}{w}$, where $w := \frac{1}{x-1} \left(= \frac{y}{1-y} = \frac{1}{z} \right)$,  such that 
\begin{align}
A_{[0]} := \begin{pmatrix} 0 & 0 & \cdots & 0 &  - s_{[0], n}  \\
1 & 0 & \cdots & 0 & - s_{[0], n-1} \\
0 & 1  & \cdots & 0 &  - s_{[0], n-2} \\
\vdots & \vdots &  \ddots & \vdots & \vdots \\ 
0 & 0 & \cdots & 1 & -s_{[0], 1} \end{pmatrix},
\hspace{5mm}
A_{[\infty]} := \begin{pmatrix} 0 & 0 & \cdots & 0 &  - s_{[\infty], n}  \\
1 & 0 & \cdots & 0 & - s_{[\infty], n-1} \\
0 & 1  & \cdots & 0 & - s_{[\infty], n-2} \\
\vdots & \vdots &  \ddots & \vdots & \vdots \\ 
0 & 0 & \cdots & 1 &  -s_{[\infty], 1} \end{pmatrix}
\end{align}
with coefficients 
$s_{[0], 1}, \cdots, s_{[0], n}$, $s_{[\infty], 1}, \cdots, s_{[\infty], n} \in k$  satisfying 
\begin{align}
\prod_{j=1}^n (\overline{\partial} - (a'_{[0]}+1 -\beta_j)) = \overline{\partial}^n + s_{[0], 1} \cdot \overline{\partial}^{n-1} + \cdots + s_{[0], n}
\end{align}
 where we set $\beta_n := 1$, and 
 \begin{align}
 \prod_{j=1}^n (\overline{\partial} - (a'_{[1]} + \alpha_j)) = \overline{\partial}^n + s_{[\infty], 1} \cdot \overline{\partial}^{n-1} + \cdots + s_{[\infty], n}.
 \end{align}
In particular,  such a log connection is uniquely determined, so
the equality $\breve{\nabla}_\circ^\diamondsuit = \breve{\nabla}_\bullet^\diamondsuit$ holds.
By  Proposition \ref{PRop532},  we conclude that  $\nabla_\circ^\diamondsuit = \nabla_\bullet^\diamondsuit$.
Thus,  $\msE^\spadesuit_\circ$ is isomorphic to $\msE^\spadesuit_\bullet$, and 
this completes the proof of this proposition.
\end{proof}

\vspace{10mm}
\section{$2$d TQFT for dormant opers in characteristic $\leq 7$} \label{S19}

In this section, we discuss the $2$-dimensional topological quantum field theory ($2$d TQFT) associated to  dormant $\mr{PGL}_n$-opers,   established in ~\cite[Theorem C, (ii)]{Wak20}.
As an application of Theorem \ref{Prop90} established  above, along with various results from prior work,
we obtain an explicit description of those TQFTs for the case   $p \leq 7$.

\LSP
\subsection{Moduli space of dormant $\mr{PGL}_n$-opers} \label{SS7}

Let $(g, r)$ be a pair of nonnegative integers with $2g-2+r >0$, 
and denote by
$\overline{\mcM}_{g, r}$ the moduli stack classifying $r$-pointed stable curves of genus $g$ over $k$.
Note that the notion  of a dormant $\mr{PGL}_n$-oper can be extended, within the framework of logarithmic geometry,  to the case where the underlying curve is a pointed stable curve (see ~\cite{Wak5} for the study of dormant $\mr{PGL}_n$-opers on such a curve).
This generalization is essential for constructing the compactified moduli stack for carrying out degeneration arguments that reduce various problems to the case of small genus.

In fact, for an $r$-tuple $\rho := (\rho_i)_{i=1}^r$ of  elements of 
$\Xi_{p, n}^{r}$
 (where $\rho := \emptyset$ if $r = 0$),
one can obtain the category
\begin{align} \label{Eq100}
\mcO p_{n, \rho, g, r}^{^\mr{Zzz...}}
\end{align}
of pairs $(\msX, \msE^\spadesuit)$ consisting of an $r$-pointed stable curve of genus $g$ over $k$  and a  dormant $\mr{PGL}_n$-oper $\msE^\spadesuit$ on $\msX$ of radii $\rho$.
According to ~\cite[Theorem C]{Wak5},  $\mcO p_{n, \rho, g, r}^{^\mr{Zzz...}}$ can be represented by a (possibly nonempty) proper Deligne-Mumford stack over $k$ and the projection 
\begin{align} \label{Eq34}
\Pi_{n, \rho, g, r} :\mcO p_{n, \rho, g, r}^{^\mr{Zzz...}} \rightarrow \overline{\mcM}_{g, r}
\end{align}
 given by $(\msX, \msE^\spadesuit) \mapsto \msX$ is finite.
Furthermore, if $p > 2n$, then it follows from ~\cite[Theorem G]{Wak5} that $\Pi_{n, \rho, g, r}$ is generically \'{e}tale, or more precisely, \'{e}tale over the points of $\overline{\mcM}_{g, r}$ classifying totally degenerate curves (cf., e.g.,  ~\cite[Definition 7.15]{Wak5} for the definition of a totally degenerate curve).
In particular, it makes sense to discuss  its generic degree
\begin{align} \label{Eq113}
N_{p, n, \rho, g, r} :=\mr{deg} (\Pi_{n, \rho, g, r}) \in \mbZ_{\geq 0}.
\end{align}
Here, 
we use the notation $(-)^\blacktriangledown$ to denote the map $\Xi_{p, n} \mapsto \Xi_{p, p-n}$ given by $A \mapsto \{ -a \, | \, a \in   \mbF_p \setminus A\}$.
The following theorem summarizes key prior results that aid in computing   the values $N_{p, n, \rho, g, r}$.

\SSP
\bt \label{Th3}
\begin{itemize}
\item[(i)]
We shall set  $\varepsilon := \LL 0, 1, \cdots, p-1 \RR$, which is the unique element of  $\Xi_{p, p-1}$.
Also, write $\rho :=  (\varepsilon, \cdots, \varepsilon) \in \Xi_{p, p-1}^{r}$.
Then,   the projection $\Pi_{p-1, \rho, g, r} : \mcO p^{^\mr{Zzz...}}_{p-1, \rho, g, r} \rightarrow \overline{\mcM}_{g, r}$ is an isomorphism.
\item[(ii)]
There exists a duality  isomorphism $\delta_{n, \rho} : \mcO p_{n, \rho, g, r}^{^\mr{Zzz...}}\xrightarrow{\sim} \mcO p_{p-n, \rho^\blacktriangledown, g, r}^{^\mr{Zzz...}}$ satisfying $\delta_{p-n, \rho^\blacktriangledown} \circ \delta_{n, \rho} = \mr{id}$.
In particular, the following equality holds
\begin{align}
N_{p, n, \rho, g, r} = N_{p, p-n, \rho^\blacktriangledown, g, r}.
\end{align}
\item[(iii)]
Suppose that $r = 0$ and $p > n \cdot \mr{max} \{ g-1, 2 \}$.
Then,  $N_{p, n, \emptyset, g, 0}$ is given  by the  formula
\begin{align}
N_{p, n, \emptyset, g, 0}  = \frac{ p^{(n-1)(g-1)-1}}{n!} \cdot  
 \sum_{\genfrac{.}{.}{0pt}{}{(\zeta_1, \cdots, \zeta_n) \in \mbC^{n} }{ \zeta_i^p=1, \ \zeta_i \neq \zeta_j (i\neq j)}}
 \frac{(\prod_{i=1}^n\zeta_i)^{(n-1)(g-1)}}{\prod_{i\neq j}(\zeta_i -\zeta_j)^{g-1}}.
\end{align}
In particular, when $n=3$, $g = 2$, and $r =0$, we have
\begin{align}
N_{p, 3, \emptyset, 2, 0} = \frac{1}{181440} \cdot p^8 + \frac{1}{4320} \cdot p^6 - \frac{11}{8640} \cdot p^4 + \frac{47}{45360}\cdot p^2.
\end{align}
\end{itemize}
\et
\begin{proof}
Assertion (i) follows from ~\cite[Theorem B]{Wak2} (see also ~\cite[Theorem A]{Hos} in the case $r = 0$).
Assertion (ii) follows from ~\cite[Theorem A]{Wak2}.
Also, 
assertion (iii) follows from  ~\cite[Theorem H]{Wak5}.
\end{proof}
\SSP

We may now reformulate  Theorem \ref{Prop90}  as the following assertion.

\SSP
\bt \label{Th6}
Let   $\rho$ be an element of $\mr{Hyp}_{p, n}$.
Then, $\mcO p^{^\mr{Zzz...}}_{n, \rho, 0, 3}$ is isomorphic to $\mr{Spec} (k)$.
In particular,
the equality $N_{p, n, \rho, 0, 3} = 1$ holds.
\et
\SSP

\LSP
\subsection{$2$d TQFT for dormant $\mr{PGL}_n$-opers} \label{SS7}

 To formulate a factorization property for  $N_{n, \rho, g, r}$'s,
 we recall the definition of a $2$d TQFT.
  For a precise and detailed account,  we refer  to ~\cite{Koc}, as well as   ~\cite{Ati}, ~\cite{DuMu1}, ~\cite{DuMu2}.

Let $\Sigma$ and $\Sigma'$ be 
closed oriented $1$-dimensional manifolds. 
 An {\bf  oriented cobordism} from $\Sigma$ to $\Sigma'$ is defined as a compact oriented $2$-dimensional manifold $M$ together with smooth maps $\Sigma\rightarrow M$, $\Sigma' \rightarrow M$
  such that $\Sigma$ maps diffeomorphically (preserving orientation) onto the in-boundary of $M$, and $\Sigma'$ maps diffeomorphically (preserving orientation) onto the out-boundary of $M$.
 We denote such a cobordism  by $M : \Sigma \Rightarrow \Sigma'$.
 Two 
 oriented cobordisms $M, M' : \Sigma \Rightarrow \Sigma'$ are {\bf equivalent} if there exists  an orientation-preserving diffeomorphism $\psi : M \xrightarrow{\sim} M'$ inducing the identity morphisms of $\Sigma$ and $\Sigma'$.
 In this way, one obtains  the category
  $2\text{-}\mcC ob$ whose  objects are 
   $1$-dimensional closed oriented manifolds and whose  morphisms 
   from $\Sigma$ to $\Sigma'$  are 
    equivalence classes
   of oriented cobordisms $M : \Sigma \Rightarrow  \Sigma'$.
   The composition of morphisms is given by gluing cobordism classes, and 
this category carries a structure of symmetric monoidal category under disjoint union.
On the other hand, 
let $\mcV ect_\mbQ$ denote the symmetric monoidal category 
   of finite-dimensional $\mbQ$-vector spaces, also with monoidal structure given by the tensor product.
 Following ~\cite[Section 1.3.32]{Koc},
a {\bf $2$-dimensional topological quantum field theory} (over  $\mbQ$), or {\bf $2$d TQFT} for short,  
is defined to be 
a symmetric monoidal functor of the form
\begin{align}
\mcZ : 2\text{-}\mcC ob \rightarrow\mcV ect_\mbQ.
\end{align}

 Note that 
 each isomorphism class of objects in $2\text{-}\mcC ob$ can be classified by an integer $n \in \mbZ_{\geq 0}$ indicating the number of connected components, i.e., the number of disjoint circles $\mbS:= \left\{ (x, y) \in \mbR^2 \, | \, x^2 + y^2 =1 \right\}$.
In other words, the full subcategory consisting of objects  $\{ \mbS^{r} \, | \, r \in \mbZ_{\geq 0} \}$, where $\mbS^{ 0} := \emptyset$ and  $\mbS^{r}$ denotes the disjoint union of $r$ copies of $\mbS$,  forms 
a skeleton of $2\text{-}\mcC ob$.
Moreover, 
each connected oriented cobordism in $2\text{-}\mcC ob$ may be represented by
 $\mbM^{r \Rightarrow  s}_g$  for some triple of nonnegative integers $(g, r, s)$,
 where $\mbM^{r \Rightarrow  s}_g$ denotes a connected, compact oriented surface with  in-boundary $\mbS^{r}$  and out-boundary $\mbS^{s}$.
According to ~\cite[Lemma 1.4.19]{Koc},  every oriented cobordism in $2\text{-}\mcC ob$ factors as a permutation cobordism, followed by a disjoint union of  $\mbM^{r \Rightarrow  s}_g$'s (for various triples  $(g, r, s)$), followed by a permutation cobordisms.
It follows  that a $2$d TQFT $\mcZ : \mcV ect_K \rightarrow 2\text{-}\mcC ob$ is uniquely determined by the $\mbQ$-vector space $A := \mcZ (\mbS^1)$ together with a collection  of $\mbQ$-linear maps 
\begin{align}
\mcZ (\mbM^{r \Rightarrow  s}_g) : A^{\otimes r}  \left(= \mcZ (\mbS^r) \right)\rightarrow  A^{\otimes s}\left(= \mcZ (\mbS^s) \right)
\end{align}
 for $(g, r, s) \in \mbZ_{\geq 0}^{3}$ (where $A^{\otimes 0} := \mbQ$).
Using this formalism  of  $2$d TQFT, we have arrived at  the following result.
 
\SSP
\bt[cf. ~\cite{Wak20}, Theorem C, (ii)] \label{Theorem4f4}
There exists a unique  $2$d TQFT
 \begin{align} \label{DeeeQQ}
 \mcZ_{n} : 2\text{-}\mcC ob  \rightarrow \mcV ect_\mbQ
 \end{align}
 determined by the following properties:
\begin{itemize}
\item
$\mcZ_{n} (\mbS^r) = (\mbQ^{\Xi_{p, n}})^{\otimes r}$, i.e., the $r$-fold tensor product of the $\mbQ$-vector space with basis $\Xi_{p, n}$;
\item
$\mcZ_{n} (\mbM_0^{0 \Rightarrow 0}) = \mr{id}_\mbQ$, and 
$\mcZ_{n} (\mbM_1^{0 \Rightarrow 0}) = \frac{(p-1)!}{n! \cdot (p-n)!} \cdot \mr{id}_\mbQ$;
\item
$\mcZ_{n} (\mbM_0^{0 \Rightarrow 1}) : \mbQ \rightarrow \mbQ^{\Xi_{p, n}}$ 
and 
$\mcZ_{n} (\mbM_0^{0 \Rightarrow 2})  : \mbQ \rightarrow    (\mbQ^{\Xi_{p, n}})^{\otimes 2}$ satisfy
\begin{align} \label{eeQwe91}
\mcZ_{n} (\mbM_0^{0 \Rightarrow 1})(1) = \varepsilon 
\hspace{5mm} \text{and} \hspace{5mm}
\mcZ_{n} (\mbM_0^{0 \Rightarrow 2}) (1) = \sum_{\lambda \in \Xi_{p, n}} \lambda \otimes \lambda^\veebar, 
\end{align}
respectively (cf. Theorem \ref{Th3}, (i), for the definition of $\varepsilon$), where $(-)^\veebar$ denotes the involution on $\Xi_{p, n}$ given by $\LL a_1, \cdots, a_n\RR \mapsto \LL -a_1, \cdots, -a_n\RR$.
\item
$\mcZ_{n} (\mbM_0^{1 \Rightarrow 0}) : \mbQ^{\Xi_{p, n}} \rightarrow \mbQ$ and 
$\mcZ_{n} (\mbM_0^{2 \Rightarrow 0}) : (\mbQ^{\Xi_{p, n}})^{\otimes 2} \rightarrow \mbQ$ satisfy
\begin{align} \label{eeQwe90}
\mcZ_{n} (\mbM_0^{1 \Rightarrow 0}) (\lambda) = 
\begin{cases} 1 & \text{if $\lambda = \varepsilon$}; \\
0 & \text{if otherwise},
\end{cases}
\hspace{5mm} \text{and} \hspace{5mm}
\mcZ_{n, \N} (\mbM_0^{2 \Rightarrow 0}) (\lambda \otimes \eta) = 
\begin{cases}
1 & \text{if $\eta = \lambda^\veebar$}; \\
0 & \text{if otherwise},
\end{cases}
\end{align}
respectively.
\item
For any triple of nonnegative integers $(g, r, s)$ with $2g-2 + r +s > 0$,
the $\mbQ$-linear map $\mcZ_{n} (\mbM_g^{r \Rightarrow s}) : (\mbQ^{\Xi_{p, n}})^{\otimes r} \rightarrow (\mbQ^{\Xi_{p, n}})^{\otimes s}$ is given  by 
\begin{align} \label{eeQQ405}
\mcZ_{n} (\mbM_g^{r \Rightarrow s}) (\bigotimes_{i=1}^r \rho_i) = \sum_{(\lambda_j)_j \in \Xi_{p, n}^{s}} N_{p, n, ((\rho_i)_i, ( \lambda_{j}^\veebar)_{j}), g, r +s} \bigotimes_{j=1}^{s} \lambda_{j}.
\end{align}
\end{itemize}
\et 
\SSP


A key feature  of this $2$d TQFT  is that it reflects the  factrization properties  of the values $N_{p, n, \rho, g, r}$, which  arises  from  
the composition of cobordism classes in $2\text{-}\mcC ob$ (cf. ~\cite[Example 6.31]{Wak20}).
For example, 
 the composition $\mbM_{g_2}^{1 \Rightarrow r_2} \circ \mbM_{g_1}^{r_1 \Rightarrow 1} = \mbM_{g_1 + g_2}^{r_1 \Rightarrow r_2}$ induces, via $\mcZ_n$, 
 the following  factorization formula:
\begin{align} \label{Eq200}
N_{p, n, (\rho_1, \rho_2), g_1 + g_2, r_1 + r_2} = \sum_{\rho_0 \in \Xi_{p, n}} N_{p, n, (\rho_1, \rho_0), g_1, r_1 +1} \cdot N_{p, n, (\rho_2, \rho_0^\veebar), g_2, r_2 +1}.
\end{align}
where $g_1$, $g_2$, $r_1$, and $r_2$ are  nonnegative integers with $2g_i -1 + r_i >0$ ($i=1,2$),
 and let $\rho_1 \in \Xi_{p, n}^{r_1}$ and $\rho_2 \in \Xi_{p, n}^{r_2}$.
Similarly,  the composition $\mbM_{g}^{2 \Rightarrow r} \circ \mbM_{0}^{0 \Rightarrow 2} = \mbM_{g+1}^{0 \Rightarrow r}$
yields another  factorization
\begin{align} \label{Eq201}
N_{p, n, \rho, g+1, r} = \sum_{\rho_0 \in \Xi_{p, n}} N_{p, n, (\rho, \rho_0, \rho_0^\veebar), g, r+2}.
\end{align}
for nonnegative integers $g$, $r$ with $2g + r >0$ and tuples $\rho \in \Xi_{p, n}^{r}$ and $\rho_0 \in \Xi_{p, n}$.
Through such recursive identifies as  \eqref{Eq200} and \eqref{Eq201}, 
the problem of determining all values $N_{p, n, \rho, g, r}$ (for  fixed $p$ and $n$) can be  reduced to computing the values in the simplest case  $(g, r) = (0, 3)$.

\LSP
\subsection{Computations of $2$d TQFTs for $p \leq 7$} \label{SS7}

In what follows, we explicitly compute $\mcZ_{n}$
  for $(3 \leq) p \leq 7$.
  As discussed  above,  this reduces to computing  $N_{p, n, \rho, 0, 3}$'s.
Let us denote by 
\begin{align}
O_{p, n}^{}
\end{align}
 the set of triples $\rho \in \Xi_{p, n}^{3}$  for which the stack $\mcO p^{^\mr{Zzz...}}_{n, \rho, 0, 3}$ (in characteristic $p$) is nonempty.
In particular, Proposition \ref{Prop4} implies that   $\mr{Hyp}_{p, n} \subseteq O_{p, n}$.
Since $\mcO p^{^\mr{Zzz...}}_{n, \rho, 0, 3}$ is \'{e}tale over $\mr{Spec}(k)$, 
the value $N_{p, n, \rho, 0, 3}$ counts  the number of (isomorphism classes of) dormant $\mr{PGL}_n$-opers on $\msP$ of radii $\rho$.

\SSP
\subsubsection{The case $(p, n) = (3, 2)$:}
By definition, the set $\Xi_{3, 2}$ consists of a single element, i.e.,  $\LL 0, 1 \RR$.
Hence,  Theorem \ref{Th3}, (i), implies  that $O_{3, 2} = \{ (\LL 0, 1 \RR, \LL 0, 1 \RR, \LL 0, 1 \RR) \}$.
The same  theorem also says that 
\begin{align}
N_{3, 2, \rho, 0, 3}  := \begin{cases}
1 & \text{if $\rho \in O_{3, 2}$}; \\
0 & \text{if otherwise}.
 \end{cases}
\end{align}


\subsubsection{The case $(p, n) = (5, 2)$:}
The set $\Xi_{5, 2}$ consists of  two elements $\{ \LL  0, 1\RR, \LL 0, 2 \RR \}$.
It follows that  any dormant $\mr{PGL}_2$-oper on $\msP$ in characteristic $5$ is of hypergeometric type.
Hence, one can apply Proposition \ref{Prop4} to  conclude that
the set $O_{5, 2}$ consists of the following $5$  triples:
\begin{align}
& (\LL 0, 1 \RR, \LL 0, 1 \RR, \LL 0, 1 \RR), \hspace{3mm}
(\LL 0, 1 \RR, \LL 0, 2 \RR, \LL 0, 2 \RR),\hspace{3mm}
(\LL 0, 2 \RR, \LL 0, 1 \RR, \LL 0, 2 \RR),\hspace{3mm}
(\LL 0, 2 \RR, \LL 0, 2 \RR, \LL 0, 1 \RR),  \\
& (\LL 0, 2 \RR, \LL 0, 2 \RR, \LL 0, 2 \RR).
\end{align}
This result also follows from the discussion in ~\cite[Section 1.6]{Iha} for $p =5$.
Moreover, Theorem \ref{Th3} (or ~\cite[Introduction, Theorem 1.3]{Moc2}, ~\cite[Theorem 10.13]{Wak20}) says  that
\begin{align}
N_{5, 2, \rho, 0, 3}  := \begin{cases}
1 & \text{if $\rho \in O_{5, 2}$}; \\
0 & \text{if otherwise}.
 \end{cases}
\end{align}


\subsubsection{The case $(p, n) = (5, 3)$:}
The set $\Xi_{5, 3}$  is given by $\{ \LL  0, 1, 2\RR, \LL 0, 1, 3 \RR\}$.
According to  Theorem \ref{Th3}, (ii), there exists a duality 
 between  dormant $\mr{PGL}_{3}$-opers and   dormant $\mr{PGL}_2$-opers.
 In particular, the set $O_{5, 3}$ can be obtained by applying the involution $(-)^\blacktriangledown$ to the elements of  $O_{5, 2}$.
 Hence,  it consists of the following $5$ triples:
\begin{align}
&(\LL 0, 1, 2\RR), \LL 0, 1, 2 \RR, \LL 0, 1, 2 \RR), \hspace{3mm}
 (\LL 0, 1, 2\RR), \LL 0, 1, 3 \RR, \LL 0, 1, 3 \RR),\hspace{3mm}
  (\LL 0, 1, 3\RR), \LL 0, 1, 2 \RR, \LL 0, 1, 3 \RR),\\
 &  (\LL 0, 1, 3\RR), \LL 0, 1, 3 \RR, \LL 0, 1, 2 \RR), \hspace{3mm}
    (\LL 0, 1, 3\RR), \LL 0, 1, 3 \RR, \LL 0, 1, 3 \RR).
\end{align}
The duality also enables us to  compute the values  $N_{5, 3, \rho, 0, 3}$ from the corresponding values for  the case $(p, n) = (5, 2)$, as follows:
\begin{align}
N_{5, 3, \rho, 0, 3}  := \begin{cases}
1 & \text{if $\rho \in O_{5, 3}$}; \\
0 & \text{if otherwise}.
 \end{cases}
\end{align}


\subsubsection{The case $(p, n) = (5, 4)$:}
The set $\Xi_{5, 4}$  contains only  a single element, i.e.,  $\LL 0, 1, 2, 3 \RR$.
By Theorem \ref{Th3}, (i),  the equality 
$O_{5, 4} = (\LL 0, 1, 2, 3 \RR, \LL 0, 1, 2, 3 \RR, \LL 0, 1, 2, 3 \RR)$ holds and 
\begin{align}
N_{5, 4, \rho, 0, 3}  := \begin{cases}
1 & \text{if $\rho \in O_{5, 4}$}; \\
0 & \text{if otherwise}.
 \end{cases}
\end{align}



\subsubsection{The case $(p, n) = (7, 2)$:}
The set $\Xi_{7, 2}$ consists of $3$ elements  $\LL 0, 1\RR$, $\LL 0, 2 \RR$, and $\LL 0, 3 \RR$.
As in the case of $p=2$,
 Proposition \ref{Prop4} (or the discussions in ~\cite{Moc2} and ~\cite{Wak20})  implies  that the set $O_{7, 2}$   contains precisely  the following $14$ triples:
\begin{align}
& (\LL 0, 1 \RR, \LL 0, 1 \RR, \LL 0, 1 \RR),\hspace{3mm}
  (\LL 0, 1 \RR, \LL 0, 2 \RR, \LL 0, 2 \RR),\hspace{3mm}
  (\LL 0, 1 \RR, \LL 0, 3 \RR, \LL 0, 3 \RR),\hspace{3mm}
  (\LL 0, 2 \RR, \LL 0, 1 \RR, \LL 0, 2 \RR), \\
&  (\LL 0, 2 \RR, \LL 0, 2 \RR, \LL 0, 1 \RR), \hspace{3mm}
  (\LL 0, 2 \RR, \LL 0, 2 \RR, \LL 0, 3 \RR),\hspace{3mm}
   (\LL 0, 2 \RR, \LL 0, 3 \RR, \LL 0, 2 \RR),\hspace{3mm}
    (\LL 0, 2 \RR, \LL 0, 3 \RR, \LL 0, 3 \RR), \\
   &  (\LL 0, 3 \RR, \LL 0, 1 \RR, \LL 0, 3 \RR), \hspace{3mm}
       (\LL 0, 3 \RR, \LL 0, 2 \RR, \LL 0, 2 \RR),\hspace{3mm}
         (\LL 0, 3 \RR, \LL 0, 2 \RR, \LL 0, 3 \RR),\hspace{3mm}
           (\LL 0, 3 \RR, \LL 0, 3 \RR, \LL 0, 1 \RR), \\
  &           (\LL 0, 3 \RR, \LL 0, 3 \RR, \LL 0, 2 \RR), \hspace{3mm}
               (\LL 0, 3 \RR, \LL 0, 3 \RR, \LL 0, 3 \RR).
\end{align}
Since any dormant $\mr{PGL}_2$-oper on $\msP$ in characteristic $7$  is of hypergeometric type,
Theorem   \ref{Th3} shows that 
\begin{align}
N_{7, 2, \rho, 0, 3}  := \begin{cases}
1 & \text{if $\rho \in O_{7, 2}$}; \\
0 & \text{if otherwise}.
 \end{cases}
\end{align}


\subsubsection{The case $(p, n) = (7, 3)$:}
Note that $\Xi_{7, 3}= \left\{ w_1, w_2, w_3, w_4, w_5\right\}$, where 
\begin{align}
w_1 := \LL 0, 1, 2 \RR, \hspace{3mm}
w_2 := \LL 0, 1,  3\RR,  \hspace{3mm}
w_3 := \LL 0, 1, 4 \RR,  \hspace{3mm}
w_4 := \LL 0, 1, 5 \RR,  \hspace{3mm}
w_5 := \LL 0, 2, 4 \RR.
\end{align}
Let   $C$ be the subset of $\Xi_{7, 3}^{3}$ 
consisting of all triples $\rho$ which arises as the radii  of dormant $\mr{PGL}_3$-opers on $\msP$  of the form $\msE^\spadesuit_{\alpha, \beta}$ for some $\alpha \in \mbF_p^{n}$ and $\beta \in \mbF_p^{n-1}$.
By applying Proposition \ref{Prop4} to the case of $(p, n) = (7, 3)$,
we obtain the following explicit list of elements in $C$:
\begin{align}
&(w_1, w_1, w_1), \hspace{3mm}
 (w_1, w_2, w_4), \hspace{3mm}
 (w_1, w_3, w_3),\hspace{3mm}
(w_1, w_4, w_2), \hspace{3mm}
(w_2, w_1, w_4),\hspace{3mm}
(w_2, w_2, w_2), \\
&(w_2, w_2, w_3), \hspace{3mm}
(w_2, w_3, w_2),\hspace{3mm}
(w_2, w_3, w_5),\hspace{3mm}
(w_2, w_4, w_1),\hspace{3mm}
(w_2, w_4, w_5),\hspace{3mm}
(w_3, w_1, w_3),  \notag \\
&
(w_3, w_2, w_2), \hspace{3mm}
(w_3, w_2, w_5),\hspace{3mm}
(w_3, w_3, w_1),\hspace{3mm}
(w_3, w_3, w_3),\hspace{3mm}
(w_3, w_3, w_5),\hspace{3mm}
(w_3, w_4, w_4),\notag \\
&(w_3, w_4, w_5), \hspace{3mm}
(w_4, w_1, w_2),\hspace{3mm}
(w_4, w_2, w_1),\hspace{3mm}
(w_4, w_2, w_5),\hspace{3mm}
(w_4, w_3, w_4),\hspace{3mm}
(w_4, w_3, w_5), \notag \\
&
(w_4, w_4, w_3), \hspace{3mm}
(w_4, w_4, w_4),\hspace{3mm}
(w_5, w_1, w_5),\hspace{3mm}
(w_5, w_2, w_3),\hspace{3mm}
(w_5, w_2, w_4),\hspace{3mm}
(w_5, w_2, w_5), \notag \\
&
(w_5, w_3, w_2), \hspace{3mm}
(w_5, w_3, w_3),\hspace{3mm}
(w_5, w_3, w_4),\hspace{3mm}
(w_5, w_3, w_5),\hspace{3mm}
(w_5, w_4, w_2),\hspace{3mm}
(w_5, w_4, w_3), \notag \\
&
(w_5, w_4, w_5).
\end{align}
Since the equality 
\begin{align}
\mr{Hyp}_{7, 3} = \left\{ (\rho_{\sigma (1)}, \rho_{\sigma (2)}, \rho_{\sigma (3)}) \, | \, (\rho_1, \rho_2, \rho_3) \in C, \sigma \in \mfS_3\right\}
\end{align}
holds,
the above list implies that  $\sharp (\mr{Hyp}_{p, n}) = 52$.
As asserted in Theorem \ref{Th6},
the equality $N_{7, 3, \rho, 0, 3} =1$ holds for  $\rho \in \mr{Hyp}_{7, 3}$.

Next, 
let us compute the values $N_{7, 2, \rho, 0, 3}$ for $\rho \in O_{p, n}$.
By Theorem \ref{Th6},  it suffices to consider the case  $\rho = (w_5, w_5, w_5)$.
To  this end, recall from 
 Theorem \ref{Th3}, (i), that
the value  $N_{p, 3, \emptyset, 2, 0}$ can be computed by 
 \begin{align}
N_{p, 3, \emptyset, 2, 0} =  \left( \frac{p^8}{181440} + \frac{p^6}{4320} - \frac{11p^4}{8640} + \frac{47 p^2}{45360} \right)\Biggl|_{p=7} = 56.
\end{align}
On the other hand, suitable  factorizations as  in \eqref{Eq200} and \eqref{Eq201} yield a decomposition of this value:
\begin{align}
 N_{7, 3, \emptyset, 2, 0} & =
 \sum_{\rho \in  O_{7, 3}} N_{7, 3, \rho, 0, 3} \cdot  N_{7, 3, \rho^\veebar, 0, 3} \\
 & =  \sum_{\rho \in  \mr{Hyp}_{7, 3}}  N_{7, 3, \rho, 0, 3} \cdot  N_{7, 3, \rho^\veebar, 0, 3}
 +  \sum_{\rho \in  O_{7, 3}\setminus \mr{Hyp}_{7, 3}}  N_{7, 3, \rho, 0, 3} \cdot  N_{7, 3, \rho^\veebar, 0, 3} \notag \\
 & = \sum_{\rho \in  \mr{Hyp}_{7, 3}}  1 \cdot 1 +  \sum_{\rho \in  O_{7, 3}\setminus \mr{Hyp}_{7, 3}}  N_{7, 3, \rho, 0, 3} \cdot  N_{7, 3, \rho^\veebar, 0, 3} \notag \\
 & = 52 +  \sum_{\rho \in  O_{7, 3}\setminus \mr{Hyp}_{7, 3}}  N_{7, 3, \rho, 0, 3} \cdot  N_{7, 3, \rho^\veebar, 0, 3}.
\end{align}
It follows  that 
$O_{7, 3}\setminus \mr{Hyp}_{7, 3}$ is nonempty and only the possible element belonging to this set is $(w_5, w_5, w_5)$.
Since $(w_5, w_5, w_5)^\veebar = (w_5, w_5, w_5)$,
the above sequence of equalities implies  $N_{7, 3, (w_5, w_5, w_5), 0, 3}^2 = 56-52 =4$, i.e., $N_{7, 3, (w_5, w_5, w_5), 0, 3} =2$.
(One of the two dormant $\mr{PGL}_3$-opers of radii  $(w_5, w_5, w_5)$ can be constructed as the second symmetric power  of the unique dormant $\mr{PGL}_2$-oper of radii $(\LL 0, 2\RR, \LL 0, 2 \RR, \LL 0, 2 \RR)$.)
Summarizing this observation, we have
\begin{align}
N_{3, \rho, 0, 3} = \begin{cases} 1 & \text{if $\rho \in O_{7, 3}$}; \\
2 & \text{if $\rho = (w_5, w_5, w_5)$};
\\
0 & \text{if otherwise}.
\end{cases}
\end{align}


\subsubsection{The case $(p, n) = (7, 4)$:}
The elements of  $\Xi_{7, 4}$ are given by 
\begin{align}
v_1 := \LL 0, 1, 2, 3\RR, \hspace{3mm}
v_2 := \LL 0, 1, 2, 4\RR, \hspace{3mm}
v_3 := \LL 0, 1, 2, 5\RR, \hspace{3mm}
v_4 := \LL 0, 1, 3, 4\RR, \hspace{3mm}
v_5 := \LL 0, 1, 3, 5\RR.
\end{align}
Note that $w_1^\blacktriangledown= v_1$, $w_2^\blacktriangledown = v_2$, $w_3^\blacktriangledown = v_4$, $w_4^\blacktriangledown = v_3$, $w_5^\blacktriangledown = v_5$.
By the duality established  in Theorem \ref{Th3}, (ii),
the set $O_{7, 4}$ is obtained from $O_{7, 3}$ by applying the involution   $(-)^\blacktriangledown$ to each element.
That is,  $O_{7, 4}$ consists of the following triples:
\begin{align}
&(v_1, v_1, v_1), \hspace{3mm}
(v_1, v_2, v_4), \hspace{3mm}
(v_1, v_4, v_2), \hspace{3mm}
(v_1, v_4, v_4), \hspace{3mm}
(v_1, v_5, v_5), \hspace{3mm}
(v_2, v_1, v_4), \hspace{3mm}
(v_2, v_2, v_2), \\
& (v_2, v_2, v_4), \hspace{3mm}
(v_2, v_3, v_5), \hspace{3mm}
(v_2, v_4, v_1), \hspace{3mm}
(v_2, v_4, v_2), \hspace{3mm}
(v_2, v_4, v_5), \hspace{3mm}
(v_2, v_5, v_3), \hspace{3mm}
(v_2, v_5, v_4), \notag \\
& (v_2, v_5, v_5), \hspace{3mm}
(v_3, v_2, v_5), \hspace{3mm}
(v_3, v_3, v_3), \hspace{3mm}
(v_3, v_3, v_4), \hspace{3mm}
(v_3, v_4, v_3), \hspace{3mm}
(v_3, v_4, v_5), \hspace{3mm}
(v_3, v_5, v_4), \notag \\ 
& (v_3, v_5, v_2), \hspace{3mm}
(v_3, v_5, v_5), \hspace{3mm}
(v_4, v_1, v_2), \hspace{3mm}
(v_4, v_1, v_4), \hspace{3mm}
(v_4, v_2, v_1), \hspace{3mm}
(v_4, v_2, v_2), \hspace{3mm}
(v_4, v_2, v_5), \notag \\ 
& (v_4, v_3, v_3), \hspace{3mm}
(v_4, v_3, v_5), \hspace{3mm}
(v_4, v_4, v_1), \hspace{3mm}
(v_4, v_4, v_4), \hspace{3mm}
(v_4, v_4, v_5), \hspace{3mm}
(v_4, v_5, v_2), \hspace{3mm}
(v_4, v_5, v_3), \notag \\ 
&(v_4, v_5, v_4), \hspace{3mm}
(v_4, v_5, v_5), \hspace{3mm}
(v_5, v_1, v_5), \hspace{3mm}
(v_5, v_2, v_3), \hspace{3mm}
(v_5, v_2, v_4), \hspace{3mm}
(v_5, v_2, v_5), \hspace{3mm}
(v_5, v_3, v_2), \notag \\ 
& (v_5, v_3, v_4), \hspace{3mm}
(v_5, v_3, v_5), \hspace{3mm}
(v_5, v_4, v_2), \hspace{3mm}
(v_5, v_4, v_3), \hspace{3mm}
(v_5, v_4, v_4), \hspace{3mm}
(v_5, v_4, v_5), \hspace{3mm}
(v_5, v_5, v_1), \notag \\ 
& (v_5, v_5, v_2), \hspace{3mm}
(v_5, v_5, v_3), \hspace{3mm}
(v_5, v_5, v_4), \hspace{3mm}
(v_5, v_5, v_5).
\end{align}
Moreover, the values $N_{7, 4, \rho, 0, 3}$ are given by 
\begin{align}
N_{7, 4, \rho, 0, 3}  := \begin{cases}
1 & \text{if $\rho \in O_{7, 2} \setminus \{ (v_5, v_5, v_5) \}$}; \\
2 & \text{if $\rho = (v_5, v_5, v_5)$}; \\
0 & \text{if otherwise}.
 \end{cases}
\end{align}


\subsubsection{The case $(p, n) = (7, 5)$:}
Note that  $\Xi_{7, 5} = \{ u_1, u_2, u_3 \}$, where
\begin{align}
u_1 := \LL 0, 1, 2, 3, 4 \RR, \hspace{10mm}
u_2 :=  \LL 0, 1, 2, 3, 5 \RR, \hspace{10mm}
u_3 :=  \LL 0, 1, 2, 4, 5 \RR. \hspace{10mm}
\end{align}
The situation  is entirely  dual to the case  $(p, n) = (7, 2)$.
Under the identities  $u_1^\blacktriangledown = \LL 0, 1\RR$, $u_2^\blacktriangledown = \LL  0, 2\RR$,  and $u_3^\blacktriangledown = \LL 0, 3\RR$, the explicit description of $O_{7, 2}$ given  above shows that the $O_{7, 5}$  consists of the following $14$ triples
\begin{align}
&(u_ 1, u_1, u_1), \hspace{3mm}
(u_1, u_2, u_2), \hspace{3mm}
(u_1, u_3, u_3), \hspace{3mm}
(u_2, u_1, u_2), \hspace{3mm}
(u_2, u_2, u_1), \hspace{3mm}
(u_2, u_2, u_3), \hspace{3mm}
(u_2, u_3, u_2),  \\
&(u_2, u_3, u_3), \hspace{3mm}
(u_3, u_1, u_3), \hspace{3mm}
(u_3, u_2, u_2), \hspace{3mm}
(u_3, u_2, u_3), \hspace{3mm}
(u_3, u_3, u_1), \hspace{3mm}
(u_3, u_3, u_2), \hspace{3mm}
(u_3, u_3, u_3), \hspace{3mm}
\end{align}
Moreover,   the equalities   $N_{7, 5, \rho, 0, 3} = N_{7, 2, \rho^\blacktriangledown, 0, 3}$  shows 
\begin{align}
N_{7, 5, \rho, 0, 3}  := \begin{cases}
1 & \text{if $\rho \in O_{7, 5}$}; \\
0 & \text{if otherwise}.
 \end{cases}
\end{align}


\subsubsection{The case $(p, n) = (7, 6)$:}
The set  $\Xi_{7, 6}$ contains  precisely a single element, i.e, $\Xi_{7, 6} := \{\LL 0, 1, 2, 3, 4, 5 \RR \}$.
According to Theorem \ref{Th3}, (i), we have 
\begin{align}
O_{7, 6} = \{ (\LL 0, 1, 2, 3, 4, 5 \RR, \LL 0, 1, 2, 3, 4, 5 \RR, \LL 0, 1, 2, 3, 4, 5 \RR) \},
\end{align}
and 
\begin{align}
N_{7, 6, \rho, 0, 3}  := \begin{cases}
1 & \text{if $\rho \in O_{7, 6}$}; \\
0 & \text{if otherwise}.
 \end{cases}
\end{align}

\LSP
\subsection*{Acknowledgements}
The author was partially supported by 
 JSPS KAKENHI Grant Number 25K06933.

\vspace{10mm}


\begin{thebibliography}{99}


\bibitem[And]{And}
Y. Andr\'{e},
Sur la conjecture des $p$-courbures de Grothendieck-Katz et un probl\`{e}me de Dwork,
 {\it Geometric aspects of Dwork theory} (2004), pp. 55-112.

\bibitem[Ati]{Ati}
M. Atiyah,
Topological quantum field theories,
{\it Inst. Hautes \'{E}tudes Sci. Publ. Math.} {\bf 68} (1988), pp. 175-186.



\bibitem[BeDr]{BeDr1}
A. Beilinson, V. Drinfeld,
Quantization of Hitchin's integrable system and Hecke eigensheaves,
available at: 
{\tt https://inspirehep.net/literature/}

\bibitem[BeHe]{BeHe}
F. Beukers, G. Heckman,
Monodromy for the hypergeometric function  $_n F_{n-1}$, 
{\it Invent. Math.} {\bf 95} (1989), pp. 325-354.



\bibitem[BeTr]{BeTr}
R. Bezrukavnikov, R. Travkin,
Quantization of Hitchin integrable system via positive characteristic,
 {\it  Pure Appl. Math. Q.} {\bf 21} (2025), pp. 635-661.


\bibitem[DuMu1]{DuMu1}
O. Dumitrescu, M. Mulase,
Edge contraction on dual ribbon graphs and $2$D TQFT,
{\it J. Algebra} {\bf 494} (2018), pp. 1-27.

\bibitem[DuMu2]{DuMu2}
O. Dumitrescu, M. Mulase,
An invitation to $2$D TQFT and quantization of Hitchin spectral curves,
{\it Banach Center Publ.}, {\bf 114}
 (2018), pp. 85-144.

\bibitem[Hos]{Hos}
Y. Hoshi, 
A note on dormant opers of rank $p-1$ in characteristic $p$,
{\it Nagoya Math. J.} {\bf 235} (2019), pp. 115-126.

\bibitem[Iha]{Iha}
Y. Ihara,
Schwarzian equations,
\textit{Jour. Fac. Sci. Univ. Tokyo} Sect IA Math. {\bf 21} (1974),  pp. 97-118.


\bibitem[JoPa]{JoPa}
K. Joshi, C. Pauly,
Hitchin-Mochizuki morphism, opers and Frobenius-destabilized vector bundles over curves,
\textit{Adv. Math.} {\bf 274} (2015),  pp. 39-75.



\bibitem[JRXY]{JRXY}
K. Joshi, S. Ramanan, E. Z. Xia, J. K. Yu,
On vector bundles destabilized by Frobenius pull-back, 
{\it Compos. Math.} {\bf 142} (2006), pp. 616-630.








\bibitem[NKa1]{Kat1}
N.  M. Katz,
Nilpotent connections and the monodromy theorem: Applications of a result of Turrittin,
\textit{Inst. Hautes Etudes Sci. Publ. Math.} {\bf 39} (1970),  pp. 175-232.


\bibitem[NKa2]{Kat2}
N. M. Katz,
Algebraic solutions of differential equations ($p$-curvature and the Hodge filtration).
\textit{Invent. Math.} {\bf 18} (1972),  pp. 1-118.

\bibitem[NKa3]{NKa3}
N. M. Katz,
A conjecture in the arithmetic theory of differential equations,
{\it Bull. Soc. Math. France,} {\bf 110}, (1982), pp. 203-239.



\bibitem[NKa4]{Kat5}
N.  M. Katz,
{\it Exponential sums and differential equations},
Ann. of Math. Stud., {\bf 124}
Princeton University Press, Princeton, NJ (1990), xii+430 pp.



\bibitem[Koc]{Koc}
J. Kock,
{\it Frobenius algebras and $2D$ topological quantum field theories}
London Math. Soc. Stud. Texts, {\bf 59}
Cambridge University Press, Cambridge, (2004), xiv+240 pp.





\bibitem[LaPa]{LaPa}
Y.  Laszlo,  C. Pauly, 
The action of the Frobenius maps on rank $2$ vector bundles in characteristic $2$, 
{\it J. Algebraic Geom.} {\bf 11} (2002), pp. 219-243.


\bibitem[LiOs]{LiOs}
F. Liu, B. Osserman,
Mochizuki's indigenous bundles and Ehrhart polynomials,
\textit{J. Algebraic Combin.} {\bf 26} (2006),  pp. 125-136.













\bibitem[Moc1]{Moc1}
S. Mochizuki,
A theory of ordinary $p$-adic curves,
{\it Publ. Res. Inst. Math. Sci.} {\bf  32} (1996), pp. 957-1152.



\bibitem[Moc2]{Moc2}
S. Mochizuki,
{\it Foundations of $p$-adic Teichm\"{u}ller theory},
AMS/IP Stud. Adv. Math., {\bf 11}
American Mathematical Society, Providence, RI; International Press, Cambridge, MA (1999), xii+529 pp.



\bibitem[Mon]{Mon}
C. Montagnon,
G\'{e}n\'{e}ralisation de la th\'{e}orie arithm\'{e}tique des D-modules \`{a} la g\'{e}om\'{e}trie logarithmique, 
{\it Ph.D. thesis, L'universit\'{e} de Rennes I} (2002).
Available at: https://tel.archives-ouvertes.fr/tel-00002545



















\bibitem[Oss]{Oss3}
B. Osserman,
Logarithmic connections with vanishing $p$-curvature,
\textit{J. Pure and Applied Algebra} {\bf 213} (2009),  pp. 1651-1664.







\bibitem[Wak1]{Wak1}
Y. Wakabayashi,
An explicit formula for the generic number of dormant indigenous bundles,
\textit{Publ. Res. Inst. Math. Sci.} {\bf 50}  (2014),  pp. 383-409.


\bibitem[Wak2]{Wak2}
Y. Wakabayashi,
Duality for dormant opers,
\textit{J. Math. Sci. Univ. Tokyo} {\bf 24}  (2017),  pp. 271-320.

\bibitem[Wak3]{Wak3}
Y. Wakabayashi,
Spin networks, Ehrhart quasi-polynomials, and combinatorics of dormant indigenous bundles,
\textit{Kyoto J. Math.} {\bf 59} (2019), pp. 649-684.







\bibitem[Wak4]{Wak5}
Y. Wakabayashi,
{\it A theory of dormant opers on pointed stable curves},
Ast\'{e}risque {\bf 432}, Soc. Math. de France, (2022), ix+296 pp.



\bibitem[Wak5]{Wak11}
Y. Wakabayashi,
Topological quantum field theory for dormant opers,
\textit{arXiv: math. AG/1709.04235v3}, (2022).



\bibitem[Wak6]{Wak7}
Y. Wakabayashi,
Cyclic \'{e}tale coverings of generic curves and ordinariness of dormant opers,
{\it J. Algebra} {\bf  623} (2023), pp. 154-192.






\bibitem[Wak7]{Wak8}
Y. Wakabayashi,
A combinatorial description of the dormant Miura transformation,
{\it Bull. Iranian Math. Soc.} {\bf  49} (2023), Paper No. 81, 26 pp.



\bibitem[Wak8]{Wak9}
Y. Wakabayashi,
Inseparable Gauss maps and dormant opers,
\textit{Math. J. Okayama Univ.} {\bf 67}, (2025), pp. 1-28.



\bibitem[Wak9]{Wak10}
Y. Wakabayashi,
The generic \'{e}taleness of the moduli space of dormant $\mfs \mfo_{2\ell}$-opers,
{\it J. Geom. Phys.} {\bf  211} (2025), Paper No. 105439, 22 pp.


\bibitem[Wak10]{Wak20}
Y. Wakabayashi,
Arithmetic liftings and 2d TQFT for dormant opers of higher level,
\textit{arXiv: math. AG/2209.08528v3}, (2025).


\end{thebibliography}
\end{document}